\documentclass[11pt]{article}

\usepackage{graphicx, color, xcolor}
\usepackage{amsmath,amsfonts,amssymb, amsthm}
\def\la{\Big\langle}
\def\ra{\Big\rangle}
\def\forall{\hbox{for all}~}
\def\L{\mathbf{L}}
\def\bfv{{\mathbf v}}

\def\bfn{{\bf n}}

\def\ve{\varepsilon}

\def\A{{\cal A}}

\def\R{{\mathbb R}}

\def\implies{\Longrightarrow}
\def\vp{\varphi}
\def\c{\centerline}
\def\E{{\cal E}}
\def\H{{\cal H}}

\def\I{{\cal I}}
\def\curl{\hbox{curl}}

\def\vs{\vskip 2em}

\def\v{\vskip 1em}

\def\C{{\cal C}}
\def\J{{\cal J}}
\def\F{{\cal F}}

\def\Tilde{\widetilde}
\def\Hat{\widehat}

\def\dint{\int\!\!\int}

\def\bega{\begin{array}}
\def\enda{\end{array}}
\def\begi{\begin{itemize}}
\def\endi{\end{itemize}}

\def\ds{\displaystyle}
\def\bel{\begin{equation}\label}
\def\eeq{\end{equation}}
\def\sqr#1#2{\vbox{\hrule height .#2pt
\hbox{\vrule width .#2pt height #1pt \kern #1pt
\vrule width .#2pt}\hrule height .#2pt }}
\def\square{\sqr74}
\def\endproof{\hphantom{MM}\hfill\llap{$\square$}\goodbreak}

\newtheorem{theorem}{Theorem}[section]

\newtheorem{example}{Example}[section]
\newtheorem{proposition}{Proposition}[section]
\newtheorem{remark}{Remark}[section]

\begin{document}

\title{\bf Results and Open Questions on\\ the Boundary Control of
 Moving Sets}
\vs

\author{Alberto Bressan
\\
\, \\
Mathematics Department, Penn State University, \,\\ University Park, PA 16802, U.S.A. 
\,\\
 e-mail: ~axb62@psu.edu} 
\maketitle  
\c{\small Dedicated to Prof.~Alexander Tolstonogov in the occasion of his 85-th birthday}
\begin{abstract}
These notes provide a survey of recent results and open problems
on the boundary control of moving sets.
Motivated by the control of an invasive biological species, we consider a class of optimization problems for a
moving set $t\mapsto \Omega(t)$,  where the goal is to minimize the area of the contaminated set $\Omega(t)$ over time, plus a cost related to the control effort. 
Here the control function is the inward normal speed, assigned along the boundary $\partial \Omega(t)$. 
We also consider problems with geographical constraints, where the invasive population is restricted within an island.
Existence and structure of eradication strategies, which entirely remove 
the invasive population in minimum 
time, is also discussed.
\end{abstract}

\section{Introduction}
\label{sec:1}
\setcounter{equation}{0}
This paper contains a survey of recent results and open problems on the boundary control of moving sets.
Its main content was presented by the author in his  lectures at the 
8-th International School-Seminar on Nonlinear Analysis and Extremal Problems, 2024.

Our models were originally motivated by control
problems for reaction-diffusion equations, of the form
\bel{CRD}
u_t~=~f(u) + \Delta u - \alpha\,u.\eeq
Here  $u=u(t,x)$ is the density of an invasive population, which
grows at rate $f(u)$, diffuses in space, and can by partly removed by
implementing a control $\alpha=\alpha(t,x)$.  
For example, $u(t,x)$ may describe the density of (disease-carrying) mosquitoes at time $t$ at the location $x$, while the control
$\alpha$ describes the amount of pesticides which are used.

Given an initial 
density 
$$u(0,x)~=~\bar u(x),$$
one seeks an optimal control $\alpha$, which minimizes a functional describing
$$\hbox{[total size of the population over time] + [cost of implementing the control].}$$

For this kind of optimization problems, 
the existence of solutions and  necessary conditions for optimality can often be
achieved by  standard analytical tools. However, one rarely finds explicit formulas for the 
optimal controls.   To gain more insight, in \cite{BCS1} a simplified model was derived,
formulated in terms of the evolution of a moving set.  

\begin{figure}[ht]
\centerline{\hbox{\includegraphics[width=6cm]{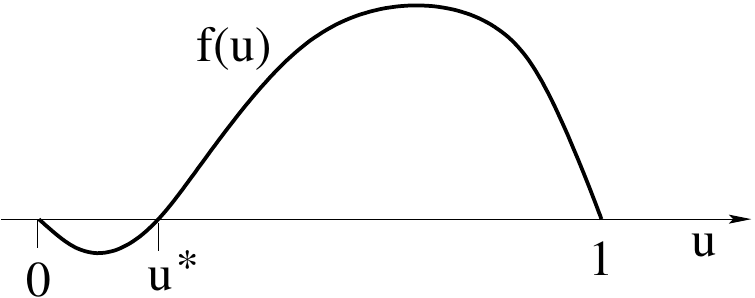}}}
\caption{\small A graph of the function $f$ in (\ref{CRD}).}
\label{f:co1}
\end{figure}

As shown in Fig.~\ref{f:co1}, assume that the  ODE
\bel{ODE1}
{d\over dt} u~=~f(u)\eeq
has two stable equilibrium points. By a variable rescaling,
we can assume these are $u=0$ and $u=1$.
We thus expect that the solution to (\ref{CRD}) will take values $u(t,x)\approx 0$ or $u(t,x)\approx 1$
on two regions, separated by a moving front.
It is thus meaningful to approximate this solution $u(t,\cdot)$ with the characteristic function
of a set 
$\Omega(t)$.    This leads to various control problems
for a moving set, which will be presently discussed.
\v

\begin{figure}[ht]
\centerline{\hbox{\includegraphics[width=13cm]{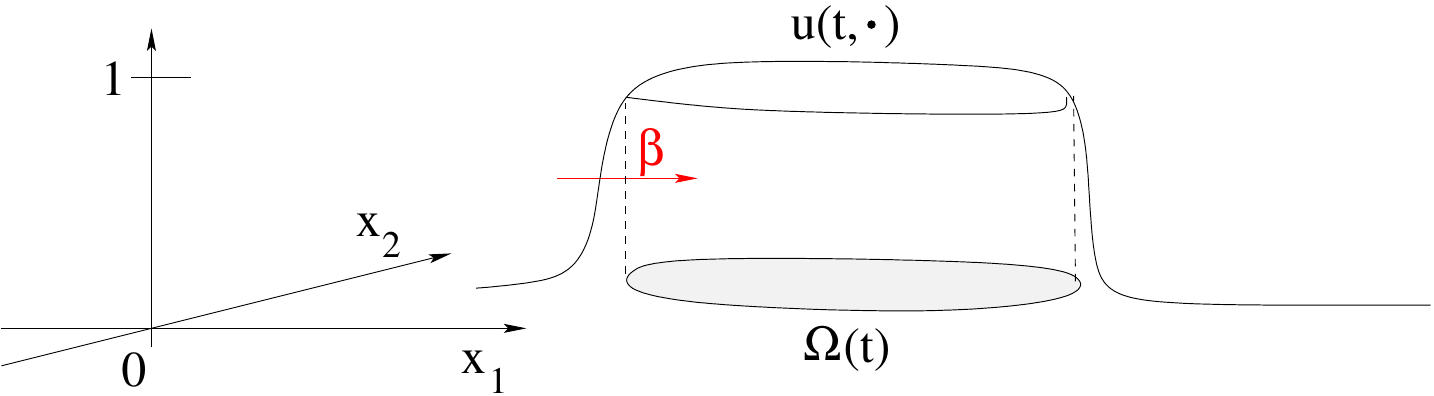}}}
\caption{\small Approximating the function $u(t,\cdot)$ with the characteristic function of a set $\Omega(t)$.
Here $\beta$ is the speed of the boundary in the inward normal direction.}
\label{f:co52}
\end{figure}

In our model, the  evolution  of a set $\Omega(t)\subset\R^2$ is controlled by
assigning the speed  $\beta(t,x) $ of each boundary point $x\in \partial \Omega(t)$ 
in the inward normal direction.
We thus regard $\beta$ as the new control variable.  

The cost for implementing 
this control will be determined
by an {\bf effort function} $E(\beta)\geq 0$.    Here $E:\R\mapsto \R_+$ is a nondecreasing
function measuring the effort needed to push boundary points inward with speed $\beta$,
per unit length of the boundary.

In an optimization problem, the goal is to reduce the size of the set $\Omega(t)$ as much as possible.
Of course, this can be done by choosing $\beta$ very large.   However,
since $E(\beta)$ is an increasing function, this control will come at a high cost.
Ultimately, the optimal strategy will be determined as a compromise between the goal of shrinking $\Omega(t)$
and avoiding a large control cost.

A crucial ingredient of our model is the choice of the effort function $E(\beta)$.   Indeed, it is this function that
provides the connection between the problem of optimal set motion and the original control problem 
for the parabolic equation (\ref{CRD}).  As explained in Section~\ref{sec:2}, this 
function can be defined in terms of an optimal control problem for traveling fronts.
\v
The remainder of the paper is organized as follows.  
In Section~\ref{sec:2} we consider a family of optimal control problems for traveling fronts.
Roughly speaking, for every $\beta\in\R$, the effort $E(\beta)$ is then defined as $\inf_\alpha \|\alpha(\cdot)\|_{\L^1}$, where $\alpha$ is a control that yields a traveling wave profile $u(t,x) = U(x-\beta t)$
with speed $\beta$, for the parabolic equation
$$u_t + f(u)~=~u_{xx} - \alpha u.$$
In Section~\ref{sec:3} we introduce optimization problems for the parabolic equation (\ref{CRD}) 
and for the set evolution problem, and show how this particular choice of the effort function 
provides a link between
the two problems.

The remaining sections deal with  the control of moving sets. 
Section~\ref{sec:4} introduces a definition 
of optimal solution, and suitable hypotheses on the effort function $E(\cdot)$, 
so that a general existence theorem
can be obtained.   Section~\ref{sec:5} contains a discussion of necessary conditions for optimality.

The last three sections cover a more specific class of problems, where the effort
function is $E(\beta)=\max\{ 0, \,\beta+1\}$.  Moreover, the  functional describing the  cost of the control
is replaced by a constraint on the total effort:
$$\E(t)~\doteq~\int_{\partial \Omega(t)} E(\beta(t,x)) \,d\sigma~\leq~ M.$$
Here the integration is w.r.t.~to the arc-length measure on the boundary of $\Omega(t)$.
Roughly speaking, this means that in absence of control the set $\Omega(t)$ expands with unit speed
in all directions. By applying a control, we can remove an area of size $M$ per unit time.
In this setting, a natural question is whether we can shrink an initial set $\Omega(0)=\Omega_0$
to the empty set in finite time, and what  strategy achieves this goal in minimum time.

The same problems can be considered in the presence of a geographical constraint. This happens when
 the invasive population
lives on an island and its expansion is restrained by a natural barrier (the sea).
In this case the moving sets
satisfy the inclusion $\Omega(t)\subseteq V$ for a fixed open set $V\subset\R^2$, and the total control effort $\E(t)$
is computed by an integral over the relative boundary $\partial \Omega(t)\cap V$.

At various places we point out the major remaining gaps in the current theory, and discuss some open
questions.

In the literature, control problems for invasive biological populations have been considered in 
\cite{ACD, ACM, ACS, CRo, SBGW, TZZ1, TZZ2}. 
In particular,  traveling waves have been studied in  \cite{HI, H1}.
For the general theory of traveling profiles for parabolic systems we refer to \cite{VVV}.
Entirely different approaches to the control of moving sets can be found in \cite{Breview,  BMN, BDZ, CLP}.

\section{Optimal control of traveling profiles}
\label{sec:2}
\setcounter{equation}{0}
In this section we consider the one dimensional parabolic equation
\bel{cpe}u_t + f(u)~=~u_{xx} - \alpha u,\eeq
under the assumptions (see Fig.~\ref{f:co1})
\begi
\item[{\bf (A1)}] {\it The function $f\in \C^2$ satisfies
\bel{f2}
f(0)~=~f(1)~=~0,\qquad f'(0)\,<\,0,\qquad f'(1)\,<\,0.\eeq
Moreover,  $f$ vanishes at only one intermediate point $u^*\in \,]0,1[\,$, where 
$f'(u^*)>0$.
}\endi
In absence of control, (\ref{cpe}) reduces to 
\bel{A1} u_t~=~f(u) + u_{xx}\,.\eeq
By definition, a traveling profile for (\ref{A1}) with speed $\beta$ is a solution of the form
\bel{TP} u(t,x)~=~U(x-\beta t).\eeq
This can be found by solving
\bel{TE} U''+\beta U' + f(U)~=~0.\eeq
Under the assumptions {\bf (A1)}, 
we seek a solution $U:\R\mapsto [0,1]$ of (\ref{TE}) with asymptotic conditions
\bel{AC}
U(-\infty) ~=~0,\qquad U(+\infty) ~=~1.\eeq
Setting $P=U'$,   a solution to (\ref{TE})-(\ref{AC}) corresponds to a heteroclinic orbit of the system
\bel{T2}
\left\{\bega{rl} U'&=~P,\\[1mm]
P'&=~-\beta P-f(U).\enda\right.\eeq
connecting the equilibrium points $(0,0)$ with $(0,1)$.
A phase plane analysis of the system (\ref{T2}) yields
\begin{theorem}\label{t:21} 
Let  $f$ satisfy {\bf (A1)}. Then
there exists a unique $\beta^*\in\R$ and a unique (up to a translation) traveling profile $U$
with speed $\beta=\beta^*$ and asymptotic values (\ref{AC}). Moreover, the function $s\mapsto U(s)$ is monotone increasing.
\end{theorem}
For a detailed proof, see Theorem~4.15 in \cite{F}.  
\v

\begin{figure}[ht]
\centerline{\hbox{\includegraphics[width=10cm]{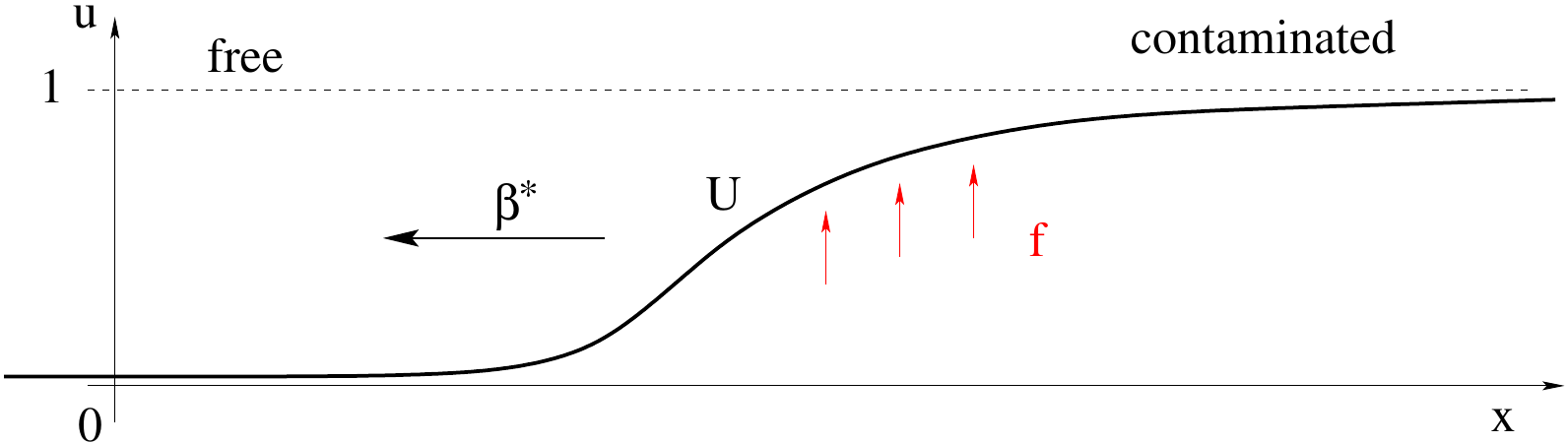}}}
\caption{\small Without any control, the parabolic equation (\ref{A1}) yields a traveling wave 
profile with speed $\beta^*<0$, and the invasive population keeps expanding.}
\label{f:co23}
\end{figure}

\begin{figure}[ht]
\centerline{\hbox{\includegraphics[width=10cm]{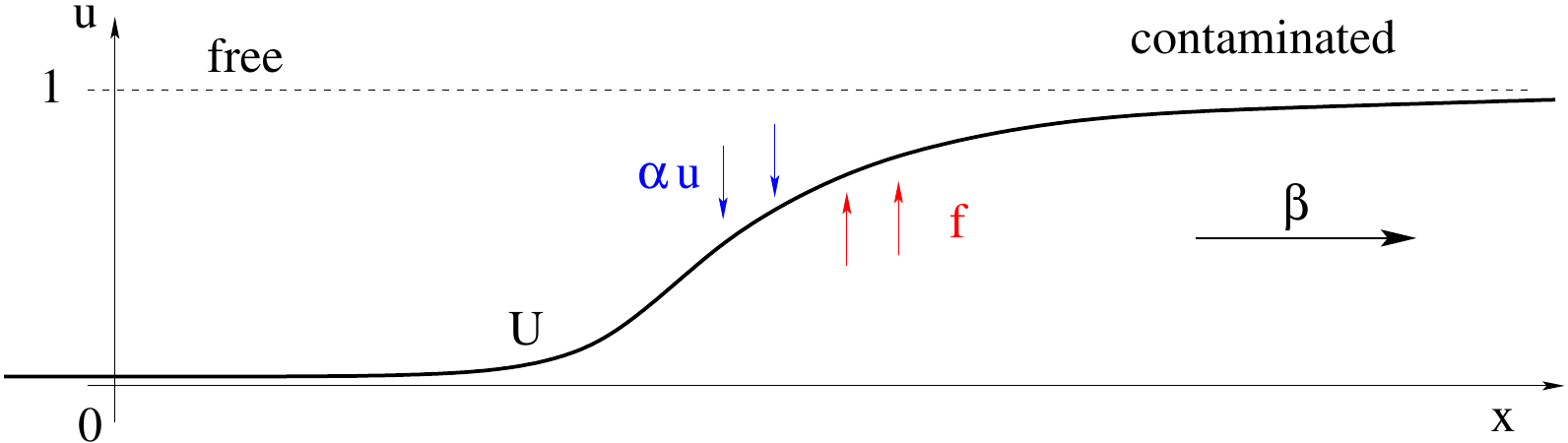}}}
\caption{\small By adding a large negative source term $\alpha u$, one can achieve a traveling profile of any 
speed $\beta >\beta^*$. In particular, when $\beta >0$ the invasive population shrinks in size.}
\label{f:co51}
\end{figure}

Since $u$ is the density of an invasive biological species, we think of the regions where  $u\approx 0$ and $u\approx 1$ as the ``free" and ``contaminated" regions, respectively.  
In absence of control, one usually has $\beta^*<0$, hence the pest 
population keeps expanding (see Fig.~\ref{f:co23}).

\begin{figure}[ht]
\centerline{\hbox{\includegraphics[width=10cm]{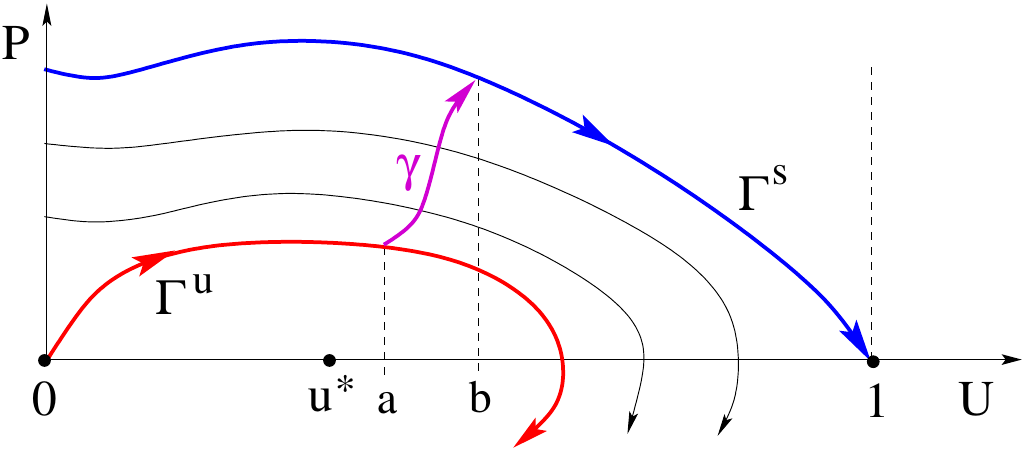}}}
\caption{\small By constructing a feedback control $\alpha(U)$ which is strictly positive on the interval 
$[a,b]$, one obtains a trajectory $\gamma$ of (\ref{T3}) which connects the unstable manifold $\Gamma^u$
through $(0,0)$ with the stable manifold $\Gamma^s$ through $(1,0)$.   This yields a traveling profile
for (\ref{cpe}) with speed $\beta$.}
\label{f:sm61}
\end{figure}

On the other hand, if a feedback control $\alpha=\alpha(u)\geq 0$ is present, this produces a negative source term on the right hand side
of (\ref{cpe}).  By a suitable choice of the control, 
one can achieve a traveling profile with any speed $\beta
\geq \beta^*$.
More precisely, a controlled traveling profile is a solution to 
\bel{ctp} U''+\beta U' + f(U) - \alpha(U)\, U~=~0.\eeq
with asymptotic conditions (\ref{AC}).  In terms of the variables $U$ and $ P=U'$, this profile corresponds to a heteroclinic
orbit of 
\bel{T3}
\left\{\bega{rl} U'&=~P,\\[1mm]
P'&=~-\beta P-f(U)+\alpha(U) \,U,\enda\right.\eeq
connecting the equilibrium points $(0,0)$ and $(0,1)$.
As shown in Fig.~\ref{f:sm61},  the phase portrait for the planar system (\ref{T2}) 
contains an unstable manifold $\Gamma^u$ through $(0,0)$ and a stable manifold $\Gamma^s$ through $(0,1)$.
For $\beta>\beta^*$ these two manifolds do not coincide.
By constructing a suitable control $\alpha(U)\geq 0$, strictly positive on a subinterval $[a,b]\subset [0,1]$, we obtain a 
trajectory $\gamma$ of (\ref{T3}) connecting $\Gamma^u$ with $\Gamma^s$.   This yields the desired traveling profile
with speed $\beta$.

Of course, there are infinitely many ways to choose this control $\alpha(U)$.  Our next goal is to construct one such control with minimum $\L^1$ norm.  
 This leads to
\begi
\item[{\bf (TWP)}] {\bf (Traveling Wave Problem).} 
{\it Given a speed $\beta\geq \beta^*$,
among all controls $\alpha\geq 0$ that yield a traveling profile with speed $\beta$,
find one that minimizes the cost
\bel{calfa}
\I(\alpha)~\doteq~\int_{-\infty}^{+\infty} \alpha\bigl(U(x)\bigr)\, dx.\eeq
}\endi

For convenience, call $\A_\beta$ the family of all Lipschitz curves 
$s\mapsto \gamma(s) = (U(s), P(s))\in \R^2$  which are the graph 
of a solution of (\ref{T3}) with $\alpha\geq 0$ and 
such that
\bel{ad1} \gamma(0)=(0,0),\qquad \gamma(\bar s) = (1,0),
\qquad
P(s)\geq 0\qquad \forall  s\in \,]0,\bar s[\,,\eeq
for some $\bar s>0$. 

Integrating along a path $\gamma\in \A_\beta$ in the $U$-$P$ plane and recalling that $P=U'$, the cost of the 
corresponding control $\alpha\geq 0$ is computed by
\bel{J11}\bega{rl} 
I(\gamma)&\ds \dot=~\int_{-\infty}^{+\infty}  \alpha\bigl(U(x)\bigr)\, dx~=~\int_0^{\bar s}
{\alpha\bigl(U(s)\bigr)\over P(s)} \, U'(s)\, ds~=~\int_0^{\bar s}  {P'(s) + \beta P(s) + f(U(s))\over U(s)} \, ds
\\[4mm]
&\ds=~
\int_0^{\bar s}  \left[\Big( {f(U(s))\over U(s) P (s)} + {\beta\over U(s)} \Big) U'(s) +{P'(s)\over U(s)}  \right]ds
~=~\int_\gamma\left[ \Big( {f(U)\over U  P } + {\beta\over U} \Big)dU+{1\over U} dP \right].
\enda\eeq
This cost must be minimized among all admissible curves $\gamma\in \A_\beta$.
Defining the vector field
\bel{vdef} \bfv~=~\left( {f(U)\over U P}+{\beta\over U}\,,~{1\over U}\right),\eeq
by (\ref{J11}) if follows
\bel{J22} I(\gamma)~=~\int_{\gamma} \bfv.\eeq

\begin{figure}[ht]
\centerline{\hbox{\includegraphics[width=7cm]{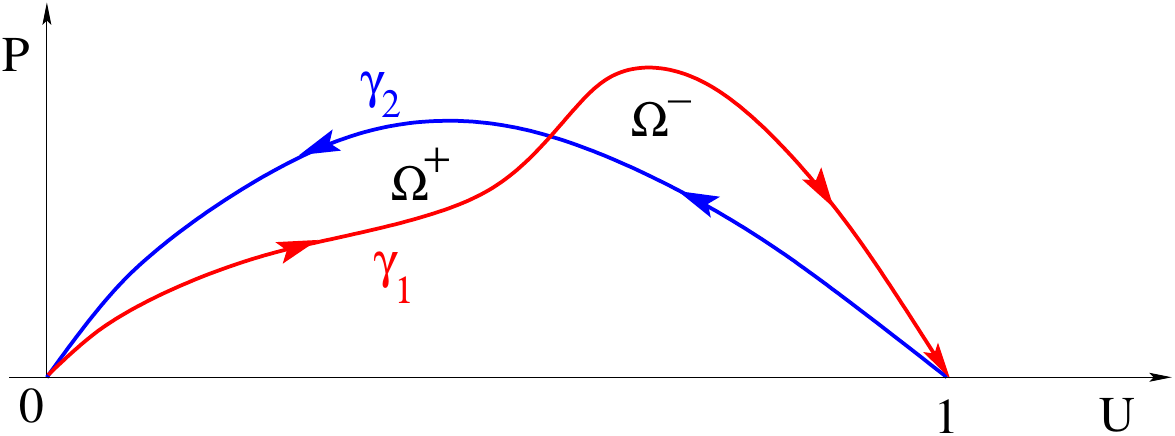}}}
\caption{\small The costs of the two paths $\gamma_1,\gamma_2$ can be compared using Stokes' theorem.}
\label{f:sm62}
\end{figure}

Following a technique developed in \cite{HH}, 
one can now use Stokes' theorem to compute the difference in cost between any two paths 
$\gamma_1,\gamma_2\in \A_\beta$. Indeed, by (\ref{J22}) we obtain
\bel{Stokes} I(\gamma_1) - I(\gamma_2)~=~
\left(\int_{\gamma_1}  - \int_{\gamma_2} \right)\bfv~=~\left( \dint_{\Omega^+} -
\dint_{\Omega^-} \right)\omega.\eeq
Here
\bel{curl2} \omega~= ~\hbox{curl}\, \bfv~=~{f(U)\over U P^2}- {1\over U^2}  \,,\eeq
while $\Omega= \Omega^+ \cup \Omega^-$ is the region enclosed between the two curves.
As shown in Fig.~\ref{f:sm62}, we call $\Omega^+$ the portion of this region whose boundary is traversed counterclockwise,
and $\Omega^-$ the portion whose boundary is traversed clockwise, when traveling first along $\gamma_1$, then along $\gamma_2$.

By (\ref{curl2}), the region where $\omega>0$ is found to be
\bel{posr2}{\cal D}^+~\doteq~
\Big\{ (U,P)\,;~\omega(U,P)>0\Big\}~=~\Big\{ (U,P)\,;~0<P< P^*(U) \Big\},\eeq
where 
\bel{P*} P^*(U)~\doteq~\sqrt{U\, f(U)}.\eeq

\begin{figure}[ht]
	\centerline{\hbox{\includegraphics[width=8.5cm]{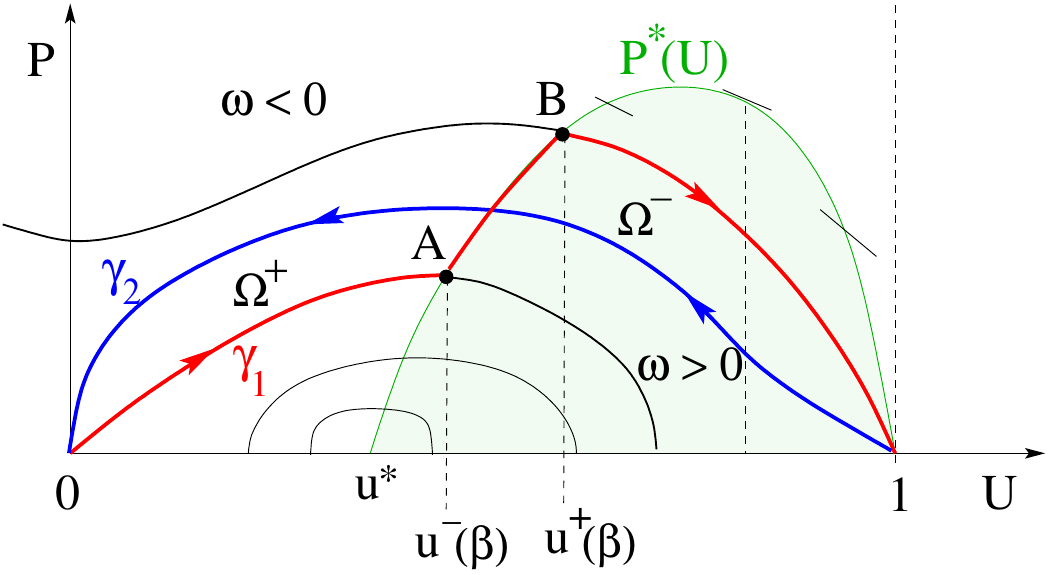}}}
	\caption{\small  By Stokes' theorem it follows that the path $\gamma_1$ going from $(0,0)$ to
		$A$, then from $A$ to $B$, then from $B$ to $(1,0)$ has a lower cost than 
		any other admissible path $\gamma_2\in \A_\beta$.  Here the control $\alpha(U)$
		is strictly positive for $U\in \bigl[u^-(\beta), u^+(\beta)\bigr]$, and vanishes elsewhere.
		Here the shaded region is the domain where $\omega<0$.}
		
	\label{f:sm63}
\end{figure}

In the situation shown in Fig.~\ref{f:sm63},  as proved in \cite{BCS1} one can
completely describe the optimal traveling profile.  Namely, the corresponding optimal path $\gamma_1$
it is a concatenation of arcs lying on three curves:
(i) the unstable manifold for the planar system (\ref{T2}) through the origin, (ii) the stable manifold for (\ref{T2}) 
through the point $(1,0)$, and (iii) the curve where $P=P^*(U)$, hence $\omega= \curl \,\bfv =0$.

\begin{theorem}\label{t:22} Let $f$ satisfy the assumptions {\bf (A1)} and let $\beta\geq \beta^*$ be given.
Let $\gamma_1$ be the path obtained by concatenating:
\begi
\item The  trajectory of (\ref{T2}) starting from $(0,0)$, until it reaches a point $A$
on the curve where $P=P^*(U)$.   
\item  The  trajectory of (\ref{T2}) through $(1,0)$,  up to the point $B$
on the curve where $P=P^*(U)$.
\item The arc of the curve where $P=P^*(U)$, between $A$ and $B$.
\endi
Assume that the above two trajectories of (\ref{T2}), passing through the points 
$A$ and $B$ respectively, do not have further
intersections with the  curve $P=P^*(U)$, for $u^*<U<1$.
Then $\gamma_1$ is optimal, i.e., it minimizes the functional $I(\gamma)$ at (\ref{J22}) among all 
admissible paths $\gamma\in \A_\beta$.
\end{theorem}
 
Indeed (see Fig.~\ref{f:sm63}), for any other admissible path $\gamma_2\in \A_\beta$, consider 
the region $\Omega = \Omega_+\cup\Omega_-$ enclosed by the concatenation of $\gamma_1 $ followed by 
$\gamma_2$. By construction, the set $\Omega^+$ is contained in the region where $\omega<0$, while
$\Omega^-$ is contained in the region where $\omega>0$.  By (\ref{Stokes}) this implies
$I(\gamma_1)\leq I(\gamma_2)$, proving the optimality of  $\gamma_1$.
\endproof

\v

Assuming that the optimal path $\gamma_1$ can be written in the form
$P=P(U)$, with $P(U)>0$ for all $0<U<1$,
a corresponding traveling wave profile $U=U(x)$ of (\ref{ctp}), in terms of the original space variable $x$,
can be obtained  by writing
\bel{xu}
x(u)~=~\int_{u^*}^u {1\over P(U)}\, dU.\eeq
Inverting the strictly increasing function $u\mapsto x(u)$ we thus obtain a traveling wave
profile
$U(x)$ satisfying (\ref{ctp}) and (\ref{AC}), where the control 
$\alpha$ achieves the minimum cost.

\begin{figure}[ht]
\centerline{\hbox{\includegraphics[width=16cm]{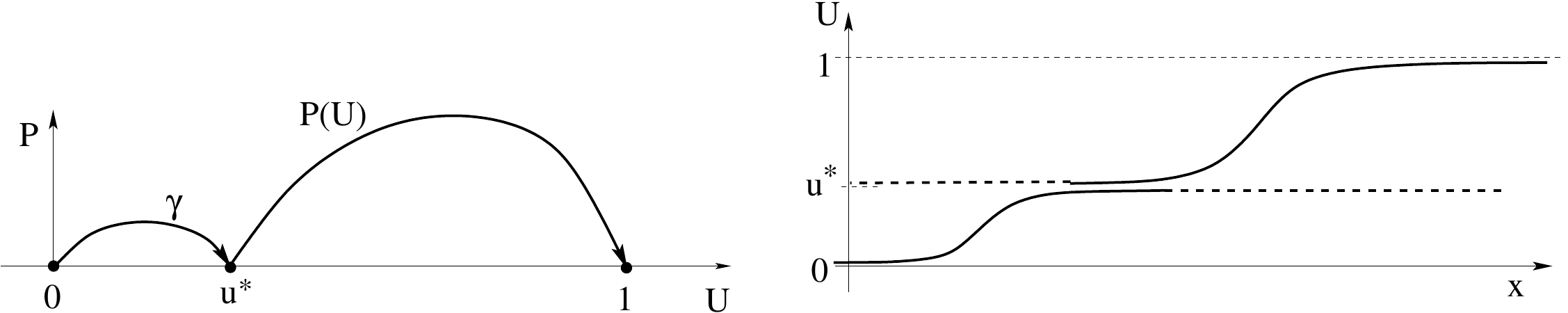}}}
\caption{\small Left: the optimal path $\gamma$, constructed  as in Theorem~\ref{t:22}, in a case where $\beta>\!> \beta^{**}$.    Right: two distinct traveling wave profiles. As explained in Remark~\ref{r:21},
in this case an optimal traveling profile for the original problem {\bf (TWP)} does not exist.}
\label{f:sm64}
\end{figure}

\begin{remark}\label{r:21} {\rm   In the case where $\beta> \beta^{**}\doteq 2\sqrt{f'(u^*)}$,
as shown in Fig.~\ref{f:sm64}, the point $(u^*,0)$ becomes a stable node for the system 
(\ref{T2}). In this case, the unstable manifold through the origin may go through the point $(u^*,0)$.
As remarked in \cite{BMZ}, if this happens the integral in (\ref{xu}) is not well defined.
One thus recovers two separate traveling wave profiles:
one joining the states
$U(-\infty)=0$ and $U(+\infty)=u^*$, and a second one joining the states
$U(-\infty)=u^*$ and $U(+\infty)=1$. These cannot be merged into a single
optimal traveling profile.  
}
\end{remark}

In all cases, the {\bf effort function} $E:[\beta^*, +\infty[\, \mapsto \R^+$  can now be defined by setting
\bel{Edef} \bega{rl} E(\beta)&\ds\doteq~\inf~\bigg\{ \int_{-\infty}^{+\infty} \alpha(x)\, dx\,;\quad 
\alpha:\R\mapsto\R_+
~~\hbox{is a control that yields}\\[1mm]
&\ds\qquad \qquad\qquad  \qquad \hbox{a traveling wave solution of (\ref{cpe}) with speed~}\beta\bigg\}.
\enda\eeq
In the situation considered in Theorem~\ref{t:22}, one has
$$E(\beta)~=~\int_\gamma \bfv,$$
where $\gamma$ is the optimal path joining $(0,0)$ with $(0,1)$ in the phase plane.
As pointed out in Remark~\ref{r:21}, however, in some cases there is no optimal control $\alpha:\R\mapsto\R_+$
with minimum $\L^1$ norm among those which yield a traveling wave with speed $\beta$.
For this reason, the infimum in (\ref{Edef}) may not always be attained as a minimum.

\v
We conclude this section by discussing the more general case of
 $n$ populations, whose spatial densities are described 
by the vector $u=(u_1, \ldots, u_n)$.
Assume that their evolution is governed by
a controlled system of reaction-diffusion equations:
\bel{psy}
u_{i,t} ~=~f_i(u,\alpha) + \sigma_i \Delta u_i\,,\qquad\qquad i=1,\ldots,n\,.\eeq
Here $\alpha=\alpha(t,x)\in\R^m$ is the control and $\sigma_1,\ldots,\sigma_n$, are diffusion coefficients.
We also assume that, in absence of control, the system of ODEs
$${d\over dt} u_i~=~f_i(u,0) ,\qquad\qquad i=1,\ldots,n\,,$$
has two asymptotically stable equilibrium points:
$u^-=(u_1^-, \ldots, u^-_n)$ and $u^+=(u_1^+, \ldots, u^+_n)$.
For a given speed $\beta\in \R$, we consider the set $\A_\beta$ of control functions
$ \alpha:\R\mapsto\R^m$ which yield a traveling wave profile with speed $\beta$, connecting the states 
$u^-, u^+$.  
That means: there exists a solution to
\bel{tv6}
f_i(U,\alpha) + \sigma_i U''_i + \beta U_i' ~=~0,\eeq
\bel{tv7} U_i(-\infty) \,=\, u_i^-,\qquad\qquad U_i(+\infty) \,=\, u_i^+,\eeq
for all $i=1,\ldots,n$.
Among all these solutions, we seek one with minimum cost, computed in terms of a Lagrangian function 
$L(u,\alpha)$. This leads to
\begi
\item {\bf Optimal Traveling Wave Problem:} 
{\it Given a speed $\beta\in \R$,  find a control $\alpha\in \A_\beta$ and 
a corresponding solution $U=U(x)$ of (\ref{tv6})-(\ref{tv7}) that  minimize the cost
\bel{calfa2}
\I(U,\alpha)~\doteq~\int_{-\infty}^{+\infty} L\bigl(U(x), \alpha(x) \bigr)\, dx.\eeq
}\endi
In the above setting, the cost functional does not admit the convenient representation
(\ref{J11}).  A fully general theory for this type of optimization problems is not yet available.
Under suitable assumptions on the functions $f$ and $L$, one can prove the existence of optimal 
solutions, and derive necessary conditions for optimality.  Some results in this direction,
valid for specific biological models, can be found in \cite{BMZ}.

\section{Controlled parabolic equations and moving sets}
\label{sec:3}
\setcounter{equation}{0}
Aim of  this section is to establish a connection between the optimal control of the parabolic 
equation (\ref{CRD}) and the control of the moving set $t\mapsto \Omega(t)$.

\begin{figure}[ht]
\centerline{\hbox{\includegraphics[width=6cm]{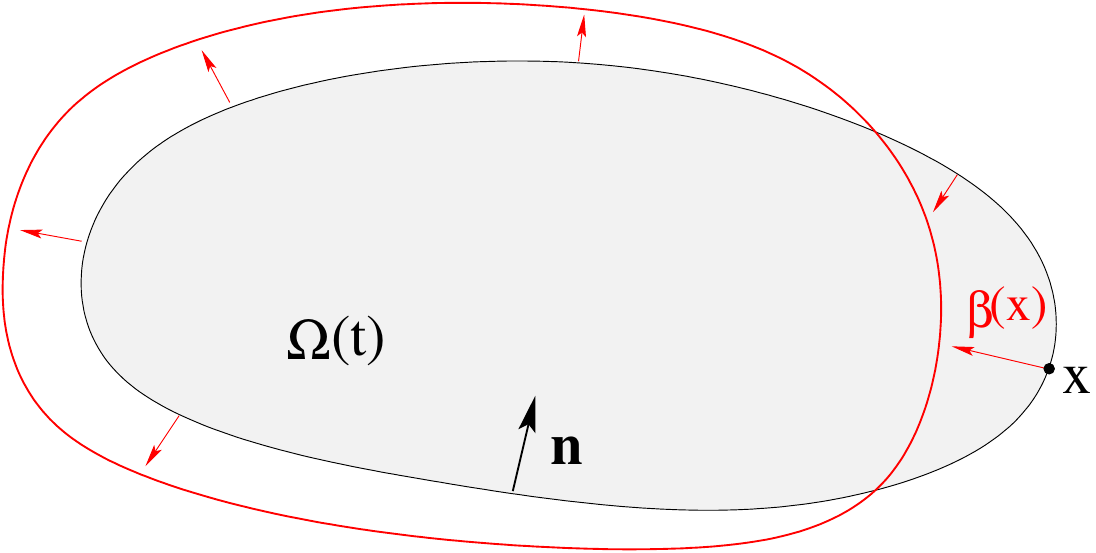}}}
\caption{\small A moving set $t\mapsto \Omega(t)$. Here the control is the function $\beta=\beta(t,x)$ determining the velocity of boundary points in the inward normal direction.}
\label{f:co27}
\end{figure}

Let $\phi:\R_+\mapsto\R_+$ be a function such that
\begi
\item[{\bf (A2)}] {\it $\phi\in \C^2$ and moreover
	\bel{pprop}
	\phi(0)\,=\,0,\qquad \phi'(0)\geq 0,\qquad \phi''(s)>0\quad\forall s>0. \eeq
}\endi
Following \cite{BCS1}, we consider two optimization problems.
\begi
\item[
{ \bf (OP1)}] {\it
Given an initial data $u(0,x)=\bar u(x)$, find a control $\alpha=\alpha(t,x)$  in (\ref{CRD})
that minimizes
\bel{OP1}\bega{c}
\ds
\J(\alpha)~\doteq~\int_0^T  \phi\left( \int  \alpha(t,x)\, dx\right) dt + \kappa_1\int_0^T \!\int u(t,x)\, dx \, dt + 
 \kappa_2\int u(T,x)\, dx \\[4mm]
\qquad =~ [\hbox{control cost}] + [\hbox{population size over time}] + [\hbox{terminal population size}].
\enda\eeq}
\endi
\begi
\item[{ \bf (OP2)}] {\it
 Given an initial set $\Omega(0)=\Omega_0$, determine a controlled evolution
$t\mapsto \Omega(t)$ which minimizes the total cost
\bel{OP2}\bega{l} J(\Omega)~\ds\doteq~\int_0^T\phi\bigl(\E(t)\bigr)\, dt + \kappa_1\int_0^T meas\bigl(\Omega(t)
\bigr)\, dt + \kappa_2\,meas\bigl(\Omega(T)\bigr)\\[2mm]
\qquad =~  [\hbox{control cost}] + [\hbox{contaminated area over time}] + [\hbox{final contaminated area}].
\enda
\eeq
}
\endi
Here the {\bf  total control effort} at time $t$ is defined as
\bel{Eff}
\E(t)~\doteq~\int_{\partial\Omega(t)}  E\bigl(\beta(t,x)\bigr)\, d\sigma,
\eeq
where $\beta(t,x)$ is the inward normal velocity at the boundary point $x\in \partial \Omega(t)$ (see Fig.~\ref{f:co27}),
$E(\beta)$ is the effort function introduced at (\ref{Edef}), and the integral 
is computed w.r.t.~to the arc-length measure along the boundary  of $\Omega(t)$.

The above two problems can be related in terms of a {\bf sharp interface limit}.
Starting with the equation
\bel{1}
u_t~=~f(u) +\Delta u - \alpha\,u\eeq
and 
rescaling the independent variables
$t\mapsto \ve t$,~ $x\mapsto \ve x$, one obtains
\bel{sil}
 u^\ve_t~=~{1\over \ve}  f(u^\ve) + \ve\, \Delta u^\ve - \alpha\,u^\ve. \eeq
Roughly speaking, the rescaling of variables amounts to looking at the spatial propagation ``from afar", 
so that the diffusion coefficient becomes smaller and 
interface between free and contaminated zone gets thinner.

The following result was proved in \cite{BCS1}. Here and in the sequel,  ${\bf 1}_S$ denotes  the characteristic function of a set $S\subset\R^2$.
\begin{theorem}\label{t:31}
For $t\in [0,T]$, let $t\mapsto \Omega(t)\subset\R^2$ denote a moving set, with $\C^1$ boundary, and inward normal speed $\beta=\beta(t,x)\geq \beta^*$.

Then there exists a family of control functions $\alpha^\ve\geq 0$ and solutions
$u^\ve$ to (\ref{sil})
such that
\bel{4}
\lim_{\ve\to 0} \Big\|u^\ve(t,\cdot) - {\bf 1}_{\Omega(t)} \Big\|_{\L^1}~=~0\quad \qquad\forall ~t\in [0,T],\eeq
\bel{5}
\lim_{\ve\to 0}  \int_{\R^2}\alpha^\ve(t,x)\, dx~=~
\int_{\partial\Omega(t)} E\bigl(\beta(t,x)\bigr)\,d\sigma.\eeq
\end{theorem}

To appreciate the meaning of this result, 
assume that $t\mapsto \Omega(t)$ is a solution to the optimization problem {\bf (OP2)}
for a moving set.  For $\ve>0$, let $u^\ve$ be a solution to the rescaled equation (\ref{sil}) and 
consider the corresponding cost functional
\bel{coep}
\J^\ve(\alpha)~\doteq~\int_0^T  \phi\left( \int  \alpha(t,x)\, dx\right) dt + \kappa_1\int_0^T \int u^\ve(t,x)\, dx \, dt + 
 \kappa_2\int u^\ve(T,x)\, dx.\eeq

Then, according to Theorem~\ref{t:31}, for $\ve>0$ small there is 
a control function $\alpha=\alpha^\ve(t,x)$ and a solution to (\ref{sil}) such that:
\begi
\item At every time $t\in [0,T]$ the density $u^\ve(t,\cdot)$ is close to the characteristic function of 
$\Omega(t)$, in the $\L^1$ distance.
\item The total cost $\J^\ve(\alpha)$ at (\ref{coep}) for the solution to the rescaled parabolic problem is almost the same as the cost $ J(\Omega)$ 
at (\ref{OP2}) for the set motion problem.
\endi

In view of these results, one may speculate if it is possible to regard the cost functional $J$
as a $\Gamma$-limit of the functionals $\J^\ve$.  Toward this goal, one would need a converse
to Theorem~\ref{t:31}.

{\bf Question 1:} {\it  In the same setting of Theorem~\ref{t:31}, consider any family of control functions  
$\alpha^\ve\geq 0$ and corresponding solutions
$u^\ve$ to (\ref{sil})
such that (\ref{4}) holds.  Does this imply that 
\bel{55}
\liminf_{\ve\to 0}  \int_{\R^2}\alpha^\ve(t,x)\, dx~\geq~
\int_{\partial\Omega(t)} E\bigl(\beta(t,x)\bigr)\,d\sigma~?\eeq
}

\begin{figure}[ht]
	\centerline{\hbox{\includegraphics[width=13cm]{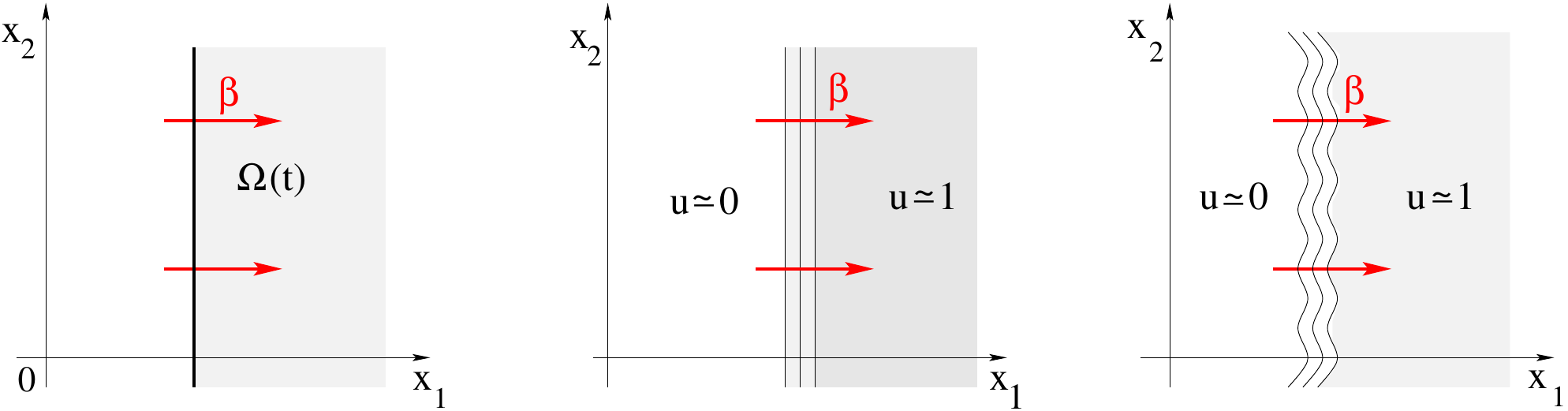}}}
	\caption{\small  Left: a set $\Omega(t)$  whose vertical boundary moves to the right, with speed $\beta$.
	Center: a traveling profile for the controlled parabolic equation, 
	moving with the same speed, constant w.r.t.~the variable $x_2$.  Right:
	a hypothetical solution of the parabolic equation, still moving to the right with the same speed $\beta$,
	but with a more complex structure, where the control  may have a lower cost.}
	\label{f:sm65}
\end{figure}

The heart of the matter is illustrated in Fig.~\ref{f:sm65}.   Assume that the boundary of the set $\Omega(t)$
is a vertical line in the $x_1$-$x_2$ plane, moving
at speed $\beta$ (Fig.~\ref{f:sm65}, left).    By implementing a control $\alpha= \alpha(x_1-\beta t)$, 
depending only on the first variable $x_1$, we can construct a traveling wave solution to (\ref{1})
which moves with the same
speed $\beta$ (Fig.~\ref{f:sm65}, center).   If we choose $\alpha$ to be the optimal control for the traveling wave problem {\bf (TWP)}, 
the cost of this control (per unit length of the boundary) is precisely $E(\beta)$.

In principle, we could produce an interface that moves with the same speed $\beta$, using a more general 
control
$\alpha=\alpha(t,x_1,x_2)$ depending on both spatial variables $x_1, x_2$  (Fig.~\ref{f:sm65}, right). 
Hence the key question is: can such a control be found, having a cost
strictly smaller than $E(\beta)$?
A negative answer to this question, although in a very special case, is provided by
Proposition~6.1 in \cite{BCS1}.

\section{Existence of optimal strategies}
\label{sec:4}
\setcounter{equation}{0}
From now on we focus on the set motion problem.
Our first goal is to reformulate the optimization problem {\bf (OP2)} in a setting 
which is general enough to achieve an existence result.
The following assumptions on 
 $E, \phi$ in (\ref{OP2})-(\ref{Eff}) will be used.
 \begi
\item[{\bf (A3)}]  {\it The function $E:\R\mapsto\R_+$ is continuous and convex.
There exist constants $\beta^*<0$ and $a>0$ such that
\bel{Eprop}
\left\{ \bega{rll}
E(\beta)&\geq ~a(\beta-\beta^*)\qquad &\hbox{if}~~\beta\geq \beta^*,\\[2mm]
E(\beta)&=~0\qquad &\hbox{if}~~\beta\leq \beta^*.\enda\right.
\eeq
In addition,  $E$ is twice continuously differentiable for $\beta>\beta^*$ and satisfies
\bel{Eass} E(\beta) - \beta\, E'(\beta)~\geq~0\qquad\qquad\forall \beta>0.\eeq
}
\endi
\begi
\item[{\bf (A4)}]  {\it The function $\phi:\R_+\mapsto\R_+
\cup\{+\infty\}$  is lower semicontinuous, nondecreasing, and convex.
Moreover, for some constants $C_1, C_2>0$ and $p>1$ one has
\bel{phia}
\phi(0)=0, \qquad \quad \phi(s)~\geq~C_1 s^p - C_2\qquad\forall ~s\geq 0\,.\eeq
}
\endi

\begin{figure}[ht]
\centerline{\hbox{\includegraphics[width=14cm]{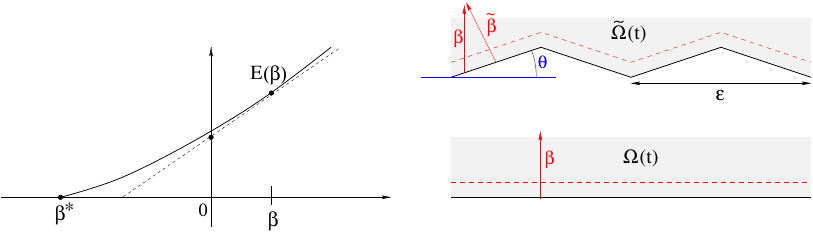}}}
\caption{\small Left: an effort function $E$ that satisfies the assumptions in {\bf (A3)}.
Right: if the assumption (\ref{Eass}) fails, the moving set $\Tilde\Omega$ with a wiggly boundary would lead to a lower cost
than the moving set $\Omega$ with smooth boundary.   As a consequence, the cost functional would not be lower
semicontinuous. 
}
\label{f:sm66}
\end{figure}

  The geometric meaning of the assumption (\ref{Eass}) is shown in Fig.~\ref{f:sm66}, left.
Namely, the tangent line to the graph of $E$ at $\beta$ has a positive intersection with the vertical axis.

This assumption is essential to guarantee the lower semicontinuity of the 
cost functional.
Indeed, consider the moving set shown in Figure~\ref{f:sm66}, right. 
Here the contaminated region $\Omega(t)$ shrinks with 
inward normal speed $\beta>0$.  We can consider an alternative strategy $\Tilde\Omega(t)$
where the boundary $\partial \Tilde\Omega(t)$ is not straight but
oscillates periodically with angle $\pm \theta$.
In the first case the cost of the control, for a portion of the boundary  parameterized by $s\in [s_1, s_2]$,  is computed by
$$\int_{s_1}^{s_2} E(\beta)\, ds,$$
where $ds$ denotes arclength along the boundary.
In the perturbed case, the boundary length increases, while the normal speed decreases.
This leads to
$$d\tilde s~ =~ {ds\over \cos\theta}\,,\qquad\qquad \Tilde \beta ~=~\beta\, \cos\theta.$$
The new total effort is
$$\int_{s_1}^{s_2} {E(\beta \cos\theta)\over \cos\theta}\, ds\,.$$
Differentiating w.r.t.~$\theta$ one obtains
$${d\over d\theta}  {E(\beta \cos\theta)\over \cos\theta}\Bigg|_{\theta=0}=~0,\qquad\qquad
{d^2\over d\theta^2}  {E(\beta \cos\theta)\over \cos\theta}\Bigg|_{\theta=0}=~E(\beta) -
\beta\, E'(\beta).$$
To achieve lower semicontinuity of the cost as $\ve\to 0$,  the wiggly profile $\partial \Tilde\Omega(t)$ should yield a larger cost. Hence
the inequality (\ref{Eass}) must be satisfied.

Another consequence of the assumptions {\bf (A3)} is that the effort function 
$E$ has sublinear growth. Indeed, when $\beta> 1$,  from (\ref{Eass}) by convexity
it follows
$$E(\beta)- \beta \,{E(\beta)-E(1)\over\beta-1}~\geq~E(\beta)- \beta \,E'(\beta)~\geq ~0.$$
Therefore
\bel{slg}E(\beta)~\leq~\beta\, E(1)\qquad\qquad \forall \beta\geq 1.\eeq
%
%
%
\v
In the following we shall
write the cost functional in a more general form, which coincides
with (\ref{OP2})-(\ref{Eff}) in the case where  the boundary $\partial \Omega(t)$ 
is smooth.
We write $\H^m$ for the $m$-dimensional Hausdorff measure,
while $\mu\Big|_V$ denotes the restriction of a measure $\mu$ to the set $V$.

\begin{figure}[ht]
\centerline{\hbox{\includegraphics[width=15cm]{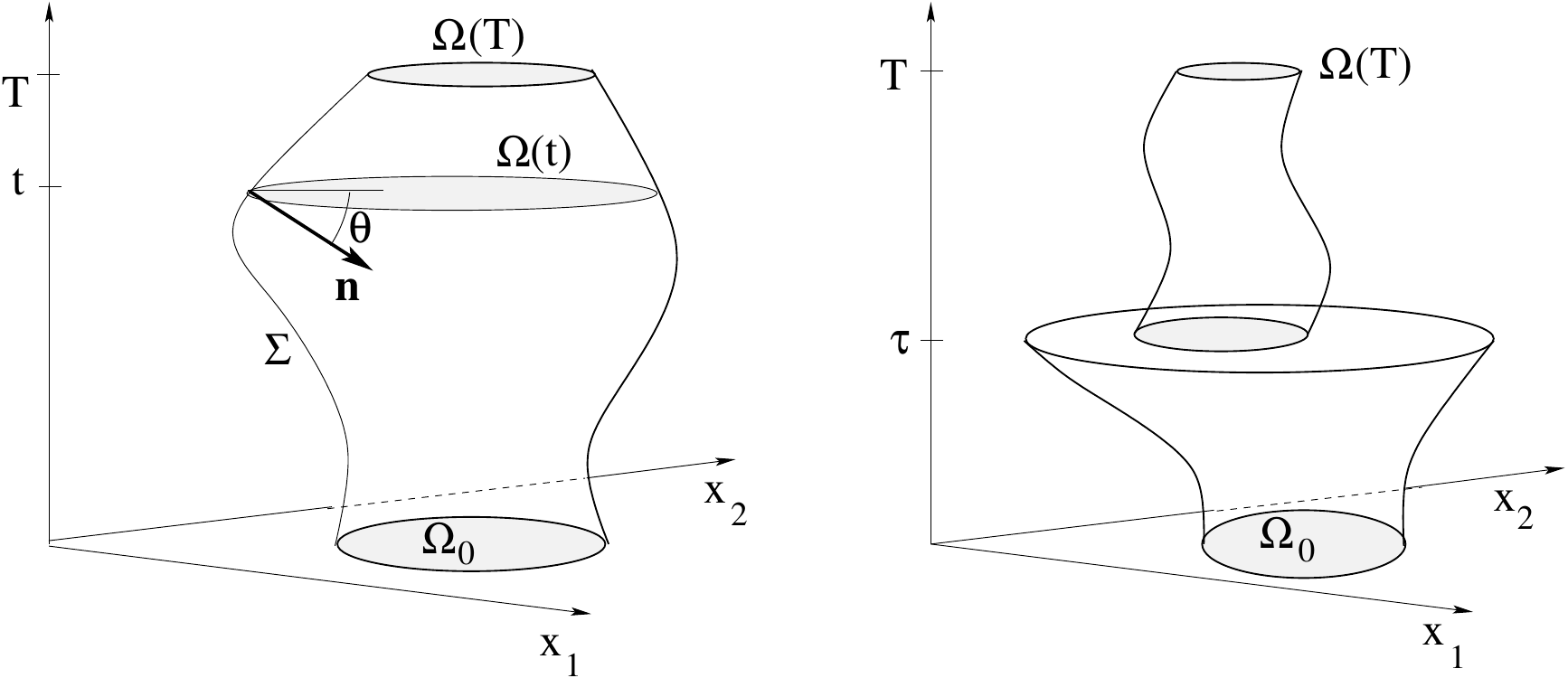}}}
\caption{\small  The measure $\mu$ defined at (\ref{mu1}) is always absolutely continuous w.r.t.~the 2-dimensional Hausdorff measure on the surface $\Sigma\subset\R^3$ where $u$ has a jump. Left: a case where
the projection $\mu^\sharp$ of $\mu$ on the $t$-axis is absolutely continuous w.r.t.~the 1-dimensional Lebesgue measure. 
Right: an example where $\mu^\sharp$ is not absolutely continuous. Indeed,
it contains a point mass at $t=\tau$.}
\label{f:sm67}
\end{figure}

Within the space $BV\bigl([0,T]\times \R^2\bigr)$ of functions with bounded variation
\cite{AFP}, 
consider the family of admissible functions 
\bel{Adef} \A~\doteq~\Big\{  
u:[0,T]\times \R^2\mapsto \{0,1\}\,;
\quad u\in BV\Big\}.\eeq  
Every $u\in\A$ determines a set 
$$\Omega~\doteq ~\bigl\{ (t,x)\in [0,T]\times \R^2\,;~u(t,x)=1\bigr\}~\subset~\R^3$$ 
with finite perimeter, and a set-valued map
\bel{sve} t~\mapsto~\Omega(t)~\doteq~\bigr\{x\in\R^2\,;~u(t,x)=1\bigl\}.\eeq
For every $u\in \A$, 
the reduced boundary $\F\Omega$  is countably ${\cal H}^2$-rectifiable.  Moreover, by Theorem~3.59 in 
\cite{AFP}, as proved by De Giorgi 
the following holds (see Fig.~\ref{f:sm67}).
\begi
\item[(i)] At each point $(t,x_1, x_2)\in \F\Omega$, the normal vector $\bfn = (n_0, n_1, n_2)\in\R^3$
is well defined.
\item[(ii)]
The distributional derivative $Du$ is a vector measure given by
\bel{Du}
Du~=~\bfn \,{\cal H}^2 \Big|_{\F\Omega}\,.\eeq
\endi
To define a cost functional for all functions $u\in\A$, 
notice that the inward normal velocity of the set $\Omega(t)$  at the boundary point 
$(t,x)\in \F\Omega$ is computed by
\bel{bdef}\beta(t,x)~=~ {-n_0\over \sqrt{n_1^2 + n_2^2}}\,.\eeq
One can thus define the scalar measure
\bel{mu1}
\mu~\doteq~ \sqrt{n_1^2+n_2^2} \cdot E\left({-n_0\over \sqrt{n_1^2 + n_2^2}}\right) \cdot \H^2\Big|_{\F\Omega}\,,\eeq
and its projection $\mu^\sharp$ on the $t$-axis, defined by
\bel{mu2}
\mu^\sharp(A)~=~\mu\Big(\{ (t,x)\,;~t\in A, ~x\in  \R^2\}\Big)\eeq
for every Borel set $A\subset [0,T]$.   We can finally define a cost functional $\Psi (u)$,
for all $u\in\A$, by setting
\bel{PU}\Psi (u)~\doteq~\left\{\bega{cl} 
\ds \int_0^T \phi\bigl(\E(t)\bigr)\, dt \qquad &  \hbox{if} ~~\mu^\sharp~\hbox{ is absolutely continuous 
with density $\E(t)$}\cr
&\qquad\qquad \qquad\qquad  \hbox{ w.r.t.~Lebesgue measure on $[0,T]$,}
\\[4mm]
+\infty 
&\hbox{if} ~~\mu^\sharp~\hbox{ is not absolutely continuous w.r.t.~Lebesgue measure.}
\enda\right.\eeq
Notice that, in  the second alternative, the infinite cost  is motivated by the assumption 
of superlinear growth in (\ref{phia}).
To justify the first alternative in (\ref{PU}) 
we observe that, in the smooth case,  the total cost of the control effort is
$$\int_0^T \phi\left(\int_{\partial\Omega(t)} E(\beta)\, d\sigma\right)dt.$$
where $d\sigma$ denotes the 1-dimensional measure along the boundary of $\Omega(t)$.
Referring to Fig.~\ref{f:sm67}, the 2-dimensional measure along the surface 
$$\Sigma~=~\bigl\{ (t,x)\,;~x\in \partial \Omega(t)\bigr\}$$  
can be expressed by
$$\H^2\Big|_\Sigma~=~d\sigma\cdot {dt\over\cos\theta}
~=~{1\over \sqrt{n_1^2+n_2^2}} \,d\sigma dt\,,
$$
where $d\sigma$ is the 1-dimensional 
measure along $\partial \Omega(t)$.
Therefore, for any $0\leq t_1 < t_2 \leq T$, we have
$$\dint_{\Sigma \cap\{ t\in [t_1, t_2]\}}  \sqrt{n_1^2+n_2^2} ~E(\beta) \,d\H^2~=~\int_{t_1}^{t_2}
 \left(\int_{\partial\Omega(t)}  E(\beta) d\sigma\right) dt\,,$$
showing that  (\ref{mu1})--(\ref{PU}) are consistent with our previous definition
of the cost for controlling the motion of the set $\Omega(t)$.

As remarked in (\ref{slg}), the effort function $E$ has sublinear growth. As a consequence,  the function 
\bel{LL} L(\bfn)~\doteq~ \sqrt{n_1^2+n_2^2} \cdot E\left( {- n_0\over \sqrt{n_1^2+n_2^2}}
\right)\eeq
is uniformly bounded.  Hence the measure $\mu$, introduced at (\ref{mu1}), is 
absolutely continuous w.r.t.~the 2-dimensional Hausdorff measure $\H^2\Big|_{\F\Omega}$ restricted to 
the set 
where $u$ has an approximate  jump.  However, the projection $\mu^\sharp$ of $ \mu$ on the time axis may 
not be absolutely continuous. For example it can have a point mass at 
a time $\tau$, as shown in 
Fig.~\ref{f:sm67}, right.  In this case, according to the  definition (\ref{PU}) the cost is $\Psi(u)=+\infty$.   

\v
In place of {\bf (OP2)}, a  generalized optimization problem for 
the moving set $t\mapsto \Omega(t)$ can now be formulated as follows.
\begi
\item[{\bf (OSM)}]  {\bf (Optimal Set Motion).} {\it Given an initial set $\Omega_0\subset\R^2$,
find an admissible function $u\in\A$, with $ u(0,\cdot)= {\bf 1}_{\Omega_0}$, 
which minimizes the functional}
\bel{FU}
\J (u)~\doteq~\Psi(u) + \kappa_1 \int_0^T\!\int u(t,x)\, dx dt + \kappa_2 \, \int u(T,x)\, dx
\,.\eeq
\endi
Here $\Psi(u)$ is the functional introduced at (\ref{PU}).  
The next theorem, on the existence of optimal solutions, was proved in \cite{BCS2}.

\begin{theorem}\label{t:41}
Let the functions $E,\phi$ satisfy the assumptions {\bf (A3)-(A4)}.
Then, for any initial set $\Omega_0\subset\R^2$ with finite perimeter and any $T>0$,
the problem {\bf (OSM)} has an optimal solution.
\end{theorem}

\section{Necessary conditions for optimality}
\label{sec:5}
\setcounter{equation}{0}
Let $t\mapsto \Omega(t)$ be an optimal solution for the problem 
{\bf (OSM)} of optimal set motion.
Aim of this section is to derive a set of necessary conditions for optimality, in a form similar to
the Pontryagin Maximum Principle.   This will require some additional assumptions on the structure and regularity
of the sets $\Omega(t)\subset\R^2$.

\begin{figure}[ht]
\centerline{\hbox{\includegraphics[width=10cm]{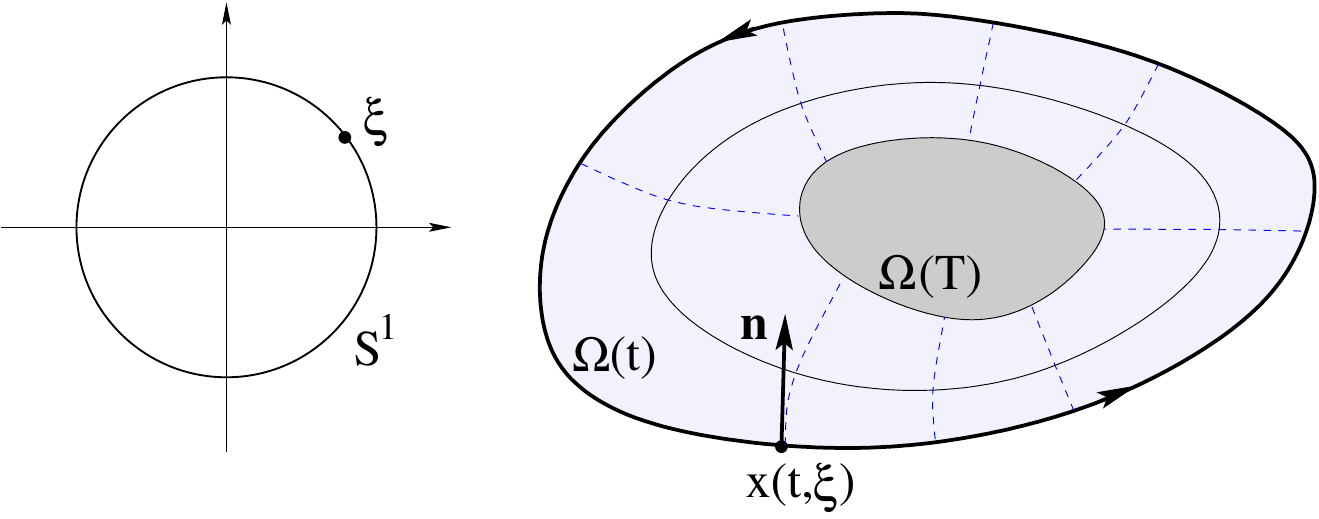}}}
\caption{\small   At each time $t\in [0,T]$, the boundary of the set $\Omega(t)$
is parameterized by $\xi\mapsto x(t,\xi)$, where $\xi\in S^1$ ranges over the unit circle.
For each $\xi$, the curve $t\mapsto x(t,\xi)$ is perpendicular to the boundaries $\partial \Omega(t)$.}
\label{f:sm70}
\end{figure}

As shown in Fig.~\ref{f:sm70},
consider the unit circumference $S^1\doteq\bigl\{\xi\in \R^2\,;~|\xi|=1\bigr\}$ and, 
for each $t\in [0,T]$, let
\bel{bpar} \xi~\mapsto ~x(t,\xi)~\in ~\partial \Omega(t)\eeq be a parameterization  of the 
boundary of $\Omega(t)$, oriented counterclockwise, satisfying the following regularity assumptions:
\begi
\item[{\bf (A5)}] {\it  The map  $(t,\xi)~\mapsto ~x(t,\xi)$ is $\C^{1,1}$ (i.e.: continuously differentiable with Lipschitz partial derivatives).
There exists a constant $C>0$ such that
\bel{ra1} {1\over C} ~\leq~\bigl|x_\xi(t,\xi)\bigr|~\leq~C\qquad\forall (t,\xi)\in [0,T]\times S^1.\eeq
Moreover, for every $\xi\in S^1$ the trajectory
$t\mapsto x(t,\xi)$ is orthogonal to the boundary $\partial\Omega(t)$ at every time $t$.   Namely,
\bel{xt1}x_t(t,\xi)~=~\beta(t,\xi)\,\bfn(t,\xi),\eeq
where $\bfn= (n_1, n_2)$ is the unit inner normal vector to $\partial \Omega(t)$ at the point $x(t,\xi)$, 
and $\beta$ is a continuous scalar function, describing the inward normal speed of the boundary.}
\endi
Throughout the following, we denote by $\bfn^\perp = (-n_2, n_1)$ the vector perpendicular to 
$\bfn=(n_1,n_2)$, and write
\bel{curv}
\omega(t,\xi)~\doteq~{1\over \bigl|x_\xi(t,\xi)\bigr|}\, \la \bfn^\perp(t,\xi),\bfn_\xi(t,\xi)\ra\eeq
for  the curvature of the boundary $\partial \Omega(t)$ at the point $x(t,\xi)$.

To derive a set of optimality conditions, we introduce the
adjoint function $Y:[0,T]\times S^1\mapsto \R$, defined as the solution of
the linearized equation
\bel{Ydt}
Y_t(t,\xi) ~=~ \left(\beta(t,\xi) -{E\bigl(\beta(t,\xi)\bigr)\over E'\bigl(\beta(t,\xi)\bigr)}\right)\omega(t,\xi) \,
Y(t,\xi)  - \kappa_1 \,,
\eeq
with terminal condition
\bel{YT} Y(T,\xi)~=~\kappa_2.\eeq
Notice that (\ref{Ydt}) yields a family of linear ODEs, that can be independently solved for
each $\xi\in S^1$.
In addition, we consider the function
\bel{lamb}\lambda(t)~\doteq~\phi'\bigl(\E(t)\bigr)~=~\phi'\left( \int_{S^1}  E(\beta(t,\xi)) \, |x_\xi(t,\xi)|\, d\xi\right).\eeq

We are now ready to state the main result, providing necessary conditions for optimality.

\begin{theorem} \label{t:51} Assume that the functions $E,\phi$ satisfy  {\bf (A3)-(A4)}. 
Let $t\mapsto \Omega(t)$ be an optimal solution to {\bf (OP2)} and assume that the boundaries 
$\partial \Omega(t)$ admit the parameterization (\ref{bpar}), satisfying {\bf (A5)}.
Then, for a.e.~$t\in [0,T]$ and $\xi\in S^1$,  the inward normal velocity $\beta=\beta(t,\xi)$ satisfies
\bel{max1}
\lambda(t) E\big(\beta(t,\xi)\bigr) -Y(t,\xi) \beta(t,\xi)  ~= ~\min_{\beta\geq \beta^*}~
\Big\{\lambda(t) E(\beta)-Y(t,\xi) \beta \Big\}.\eeq
\end{theorem}

The above result was proved in \cite{BCS2} under the stronger assumption that the map $(t,\xi)\mapsto x(t,\xi)$
is $\C^2$.  With some additional work, the same conclusion can be proved also in the $\C^{1,1}$ case.

At this stage, we must point out that a major gap still remains  in the theory.   Indeed, Theorem~\ref{t:41} 
yields the existence of optimal solutions $t\mapsto \Omega(t)$ within a general class of sets with finite perimeter.  
On the other hand, the necessary conditions require a much stronger regularity assumption.
Namely, a parameterization of the boundaries $\partial \Omega(t)$ should exist, satisfying the conditions
in {\bf (A5)}.
Notice that the $\C^{1,1}$ regularity is essential for
the orthogonal curves
$t\mapsto x(t,\xi)$ to be well defined in terms of the ODE (\ref{xt1}).  We thus state our next open problem:
\v
{\bf Problem 2.} {\it Establish the regularity of the set valued map $t\mapsto \Omega(t)$ which 
provides a solution to {\bf (OSM)}. In particular, prove that the  orthogonal curves to the boundaries 
$\partial \Omega(t)$ are well defined.
}
\v

\begin{figure}[ht]
\centerline{\hbox{\includegraphics[width=10cm]{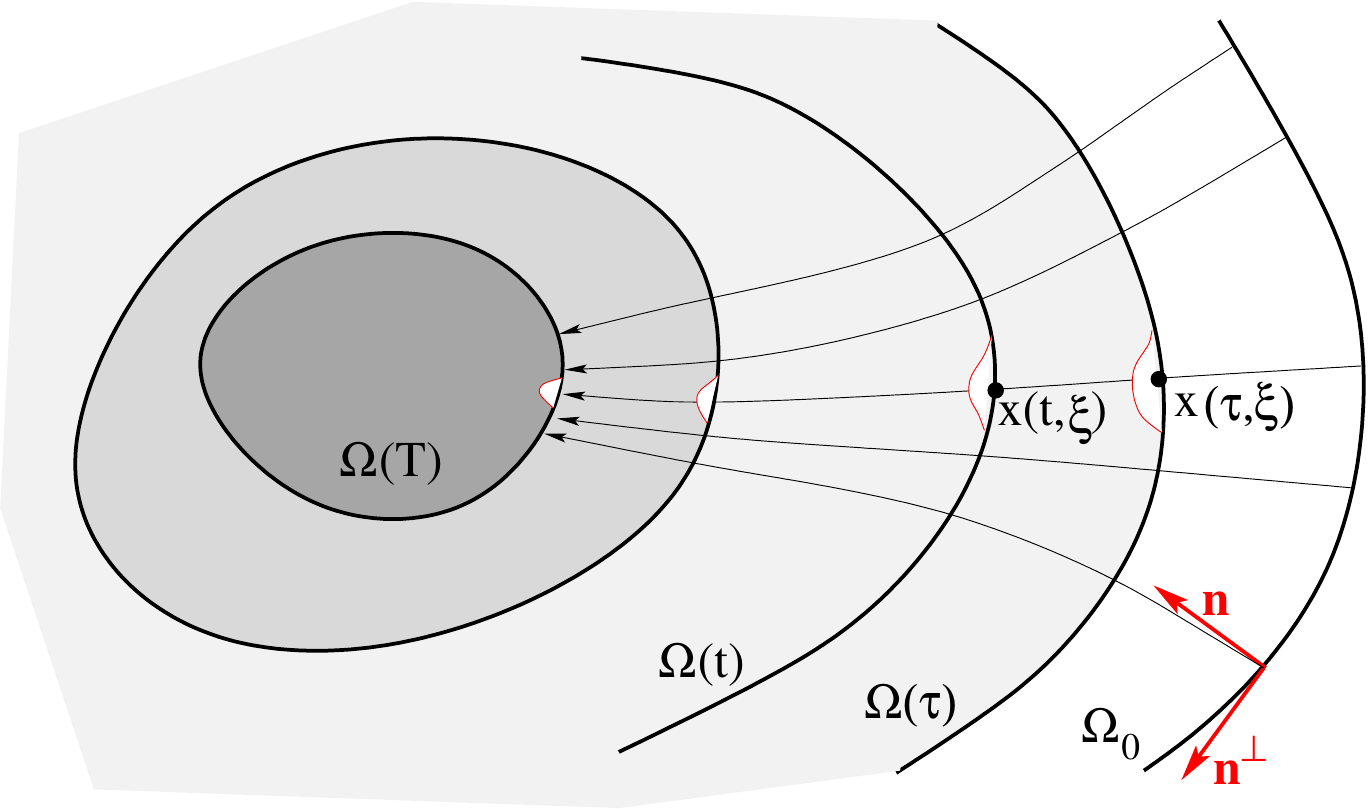}}}
\caption{\small   A localized variation of an optimal strategy.   At time $\tau$, an increase of the control effort in a neighborhood of the point $x(\tau,\xi)$ produces a  reduction of the set $\Omega(\tau)$.   If the control effort is 
kept unchanged for $t>\tau$, the corresponding sets $\Omega(t)$, $t\in [\tau, T]$ are all slightly smaller then before. }
\label{f:sm42}
\end{figure}

In the remainder of this section we explain the main ideas behind Theorem~\ref{t:51}.
As shown in Fig.~\ref{f:sm42}, assume that at some time $\tau\in [0,T]$
we change our control $\beta(\cdot)$ in a neighborhood of a point $(\tau, \xi)$.  For example, we can
locally  increase the
inward normal speed of the boundary,
thus obtaining a set $\Omega^\ve(\tau)$ which is slightly smaller than 
$\Omega(\tau)$. 
This corresponds to a ``needle variation" used in the proof of the Pontryagin Maximum Principle
\cite{BPi, FR}.
However, in the present case the perturbation is localized in time and also in space. 
Afterwards, for  $t\in [\tau, T]$, we keep the same control effort
as before at all points of the boundary $ \partial \Omega(t)$.
 As a consequence, the 
sets $\Omega(t) $ will be replaced by (possibly smaller) sets $\Omega^\ve(t)$.

By studying
the perturbed sets $\Omega^\ve(t)$ for $t\in [\tau,T]$, we can compare
\begi
\item[(1)] the additional cost of the control, needed to shrink the set $\Omega(\tau)$
in a neighborhood of 
$x(\tau, \xi)$, at time $\tau$,
\item[(2)] the reduction in the total cost, due to the fact that the measure of the sets $\Omega^\ve(t)$,
$t\in [\tau, T]$ is now smaller  than in the original solution.
\endi
If (1) is smaller than (2),  then our
perturbation will strictly decrease the sum $J(\Omega)$  of all costs in (\ref{OP2}).  This would rule out optimality.

\begin{figure}[ht]
\centerline{\hbox{\includegraphics[width=10cm]{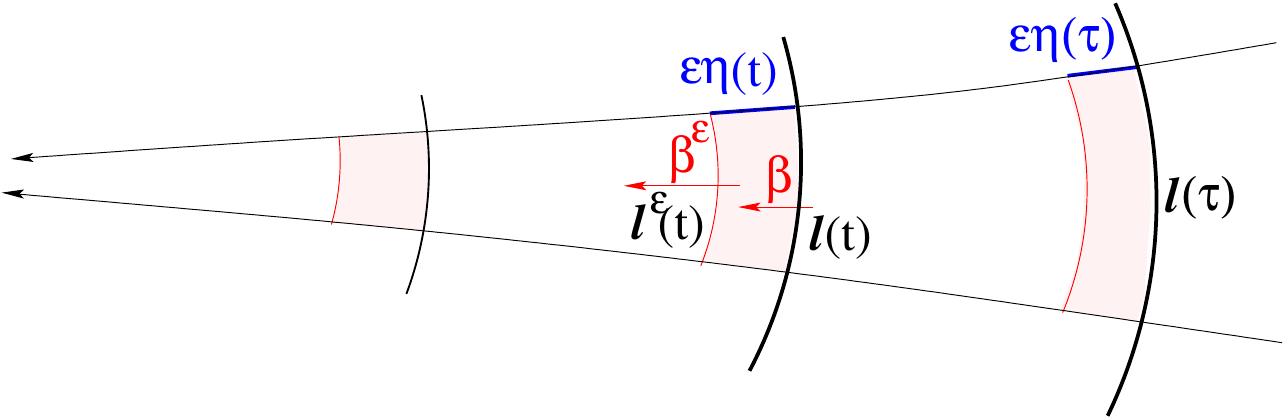}}}
\caption{\small  A tube $\Gamma$ where the sets $\Omega(t)$ are modified. 
If the curvature $\omega$ of the boundary $\partial \Omega(t)$ is positive, as  this boundary is pushed inward the
cross section decreases: $l^\ve(t)<l(t)$. As the integral of the control effort $E(\beta)$ across the tube is kept constant,
the inward speed increases: $\beta^\ve(t)>\beta(t)$.  As a consequence, the additional
infinitesimal  area  $\ve\eta(t)\, l(t)$ cleaned up in the perturbed solution increases in time.}
\label{f:sm43}
\end{figure}

To estimate the quantities in (1) and (2),  consider a small tube around the orthogonal trajectory
$t\mapsto x(\xi,t)$, say
$$\Gamma~\doteq~\bigl\{ x(t,\zeta)\,;~~t\in [\tau, T],~~\zeta\in [\xi-\delta,\xi+\delta]\bigr\}.$$
As shown in Fig.~\ref{f:sm43}, call $l(t)$ the width of the tube at time $t>\tau$, and
let $\ve \eta(t)$ the amount by which the boundary of $\Omega(t)$ is pushed inward,
in the perturbed strategy.
Notice that, if the boundary curvature $\omega$ does not vanish, to leading order we have
\bel{lep1}l^\ve(t)~=~l(t)\cdot \bigl[1-\ve \eta(t)\omega(t,\xi)\bigr].\eeq

Call $\beta^\ve$ the inward normal speed of boundary points in the perturbed strategy.
The assumption that the effort (and hence the total cost of the control) at each 
time $t>\tau$ remains the same as in the original 
motion implies
\bel{ef=1}E\bigl(\beta^\ve (t,\xi)\bigr) \cdot l^\ve(t)~=~E\bigl(\beta (t,\xi)\bigr) \cdot l(t).\eeq
Writing the first order expansion
$\beta^\ve~=~\beta + \ve b+o(\ve)$, for some leading order perturbation $b$,  this yields
$$\Big[ E\bigl(\beta (t,\xi)\bigr) + \ve E'\bigl(\beta (t,\xi)\bigr)\, b(t,\xi)\Big]
\cdot l(t)  \bigl[1-\ve \eta(t) \omega(t,\xi)\bigr]
~=~E\bigl(\beta (t,\xi)\bigr) \cdot l(t).$$
Collecting first order terms one obtains
\bel{ef=}E'\bigl(\beta (t,\xi)\bigr)\, b(t,\xi)~=~E\bigl(\beta (t,\xi)\bigr)\,\eta(t)\, \omega(t,\xi),\eeq
and hence
\bel{deta}{d\over dt} \eta(t)~=~b(t)~=~{E\bigl(\beta (t,\xi)\bigr)\,\omega(t,\xi)\over 
E'\bigl(\beta (t,\xi)\bigr)}\,\eta(t).\eeq

Next, calling $A(t)=\eta(t) l(t)$ the local change in the area,  we obtain
\bel{darea}\bega{rl}\ds {d\over dt}A(t)&\ds=~{d\over dt} \eta(t)\cdot  l(t)+\eta(t)\cdot {d\over dt} l(t)\\[4mm]
&\ds=~
{E\bigl(\beta(t,\xi)\bigr)\omega(t,\xi)\over E'\bigl(\beta(t,\xi)\bigr)}\eta(t) l(t) - \eta(t) l(t) \beta(t,\xi) \omega(t,\xi)
\\[4mm]
&\ds=~\left[{E\bigl(\beta(t,\xi)\bigr)\over E'\bigl(\beta(t,\xi)\bigr)} -\beta(t,\xi)\right]  \omega(t,\xi)A(t).\enda\eeq
Motivated by (\ref{darea}), we denote by $t\mapsto A(t,\tau,\xi)$ the solution to the ODE
\bel{Adt} {d \over d t}  A(t)~=~\left[{E\bigl(\beta(t,\xi)\bigr)\over E'\bigl(\beta(t,\xi)\bigr)} -\beta(t,\xi)\right]  \omega(t,\xi)A(t), \qquad\qquad A(\tau)=1.\eeq
Moreover, we define the adjoint function $Y=Y(t,\xi)$ to be the solution to (\ref{Ydt}) 
with terminal condition (\ref{YT}).
We claim that
\bel{Ydef} Y(\tau,\xi)~=~\int_\tau^T \kappa_1 A(t,\tau,\xi)\, dt + \kappa_2 A(T,\tau,\xi)\eeq
is the decrease in the cost, per unit area cleaned up
at time $\tau$ in a neighborhood of the boundary point $x(\tau,\xi)$.

Indeed, by (\ref{Adt}) and (\ref{Ydt}) it follows
$${\partial \over \partial t} \Big(A(t,\tau,\xi) Y(t,\xi)\Big)~=~- \kappa_1 A(t,\tau,\xi).$$
$$\int_\tau^T {\partial \over \partial t} \Big(A(t,\tau,\xi) Y(t,\xi)\Big)\, dt~=~- \int_\tau^T \kappa_1 A(t,\tau,\xi)\, dt
~=~A(T,\tau,\xi) Y(T,\xi) - A(\tau,\tau,\xi) Y(\tau,\xi).$$ 
Recalling that $A(\tau,\tau,\xi)=1$, one obtains
$$Y(\tau,\xi)~=~\int_\tau^T \kappa_1 A(t,\tau,\xi)\, dt + A(T,\tau,\xi) Y(T,\xi),$$
proving (\ref{Ydef}).

Next, consider the function
\bel{lamb2}\lambda(t)~\doteq~\phi'\bigl(\E(t)\bigr)~=~\phi'\left( \int_{S^1}  E(\beta(t,\xi)) \, |x_\xi(t,\xi)|\, d\xi\right).\eeq
Assume that the optimal control $\beta$ is continuous at the point $(\tau,\xi)$.
Consider a ``needle variation" replacing $\beta(t,\zeta)$ with some other value $\Hat\beta$
for $(t,\zeta)\in [\tau-\ve,\tau]\times [\xi-\delta, \xi+\delta]$.
To leading order, the change in the control cost will be
\bel{ccost}\lambda(\tau) \bigl[ E(\Hat\beta) - E(\beta(\tau,\xi))\bigr]\cdot \ve l(\tau).\eeq
On the other hand, the reduction in total cost (due to the fact that the sizes of all 
sets $\Omega^\ve(t)$, $t\in [\tau,T]$ become smaller) is computed by
\bel{csave}\bigl[\Hat\beta - \beta(\tau,\xi)\bigr]\cdot \ve l(\tau) \cdot Y(\tau, \xi).\eeq
If the quantity in (\ref{ccost}) were strictly smaller than (\ref{csave}), this would 
contradict optimality.  
Therefore
$$
\lambda(\tau) E\big(\beta(\tau,\xi)\bigr) -Y(\tau,\xi) \beta(t,\xi)  ~\leq~
\lambda(\tau) E(\Hat \beta)-Y(\tau,\xi) \Hat \beta.$$
Since $\Hat\beta\geq \beta^*$ is arbitrary, this yields (\ref{max1}).

\begin{remark}\label{r:51} {\rm 
The adjoint variable $Y>0$ introduced at 
(\ref{Ydt})-(\ref{YT}) can be interpreted as a ``shadow price".
Namely (see Fig.~\ref{f:sm42}),  thinking of $\Omega(t)$ as the contaminated set, 
assume that  at time $\tau$ an external contractor offered to ``clean up"
a neighborhood of the point $x(\tau, \xi)$, thus replacing the set $\Omega(\tau)$ with a smaller set 
$\Omega^\ve(\tau)$, 
 at a price of $Y(\tau,\xi)$ per unit area.   
In this case, accepting or refusing the offer would make no difference in the total cost.   
}
\end{remark}

\section{Eradication problems with geographical constraints}
\label{sec:6}
\setcounter{equation}{0}
In the remaining sections we focus on the basic case where
\bel{EP}
E(\beta)~=~\left\{ \bega{cl} 
0\quad &\hbox{if}\quad \beta\leq-1,\cr
1+\beta \quad &\hbox{if}\quad \beta> -1,\enda\right.  
\qquad\qquad \phi(s)~=~\left\{ \bega{cl} 0\quad &\hbox{if}\quad s\leq M,\cr
+\infty\quad &\hbox{if}\quad s> M.\enda\right.\eeq
Notice that now the minimization of the control cost in (\ref{OP2}) 
is replaced the constraint $\E(t)\leq M$ for a.e.~$t\in [0,T]$.
This models a situation where:
\begi
\item If the control effort is everywhere zero: $E(\beta)=0$, then the inward normal speed is $\beta= \beta^*= -1$
at every point. 
Hence the contaminated set $\Omega(t)$ expands with unit speed in all directions.
In particular, its area increases at a rate equal to the perimeter:
$${d\over dt} m_2\bigl(\Omega(t)\bigr)~=~\int_{\partial \Omega(t)} - \beta^* \, d\sigma~=~
m_1\bigl(\partial \Omega(t)\bigr).$$
\item By implementing a control with total effort  $\E(t)= M$, we can clean up a region of area $M$ per unit time:
\bel{darea2} {d\over dt} m_2\bigl(\Omega(t)\bigr)~=~\int_{\partial \Omega(t)} - \beta(t,x) \, d\sigma~=~
\int_{\partial \Omega(t)} \Big[1 - E\bigl(\beta(t,x)\bigr)\Big] \, d\sigma~=~m_1\bigl(\partial \Omega(t)\bigr)-M.\eeq
\endi
Here and in the sequel, $m_1, m_2$ denote the 1-dimensional and 2-dimensional measures. 
In this setting, two main problems will be considered.

\begi
\item[{\bf (EP)}] {\bf Eradication Problem.}  
{\it   Let an initial set $\Omega_0$ and a constant $M>0$ be given. 
Find a set-valued function $t\mapsto \Omega(t)$  such that, for some $T>0$,
\bel{nco}\Omega(0)~=~\Omega_0,\qquad\qquad 
\Omega(T)~=~\emptyset,\eeq
\bel{eb}
\int_{\partial\Omega(t)}E\bigl(\beta(t,x)\bigr)\, d\sigma~\leq~ M,\qquad \forall t\in [0,T].\eeq}
\v
\item[{\bf (MTP)}]  {\bf Minimum Time Problem}  {\it  Among all strategies that satisfy
(\ref{nco})-(\ref{eb}), find one the minimizes the time $T$.}
\endi
Thinking of $\Omega(t)$ as the region contaminated by a pest population at time $t$, we thus seek an admissible strategy that eradicates the infestation within finite (or minimum) time.

The next result provides a simple condition for the solvability of the Eradication Problem. 
A proof can be found in \cite{BBC}.

\begin{proposition} 
Assume that the the convex closure of $\Omega_0\subset\R^2$ has perimeter
\bel{pco} m_1\bigl(\partial (co \,\Omega_0)\bigr)\,<\,M.\eeq
Then the eradication problem {\bf (EP)} has a solution.
\end{proposition}
\v

The above  problems become more interesting in connection with an additional
spatial constraint.   Such a constraint arises naturally when 
the invasive population spreads over an island   $V\subset\R^2$.  In this case, a geographical barrier (the sea)
prevents any further expansion.

To model this situation, let $V\subset\R^2$ be a bounded open set.
Assuming that the entire domain $V$ is initially contaminated, we consider
\begi
\item[{\bf (CEP)}] {\bf Constrained Eradication Problem.}  
{\it   Let $M>0$ be given. 
Find a set-valued function $t\mapsto \Omega(t)\subseteq V$  such that, for some $T>0$,
\bel{era}\Omega(0)~=~V,\qquad\qquad 
\Omega(T)~=~\emptyset,\eeq
\bel{ebo}
\int_{\partial\Omega(t)\cap V}E\bigl(\beta(t,x)\bigr)\, d\sigma~\leq~ M,\qquad \forall t\in [0,T].\eeq}
\endi
Notice that in (\ref{ebo}) the integral of the effort ranges only over the relative boundary of $\partial\Omega(t)$,
inside the open set $V$.
If an eradication strategy exists, we can then seek an optimal one, minimizing the time $T$.

It is clear that the existence (or non-existence) of an eradication strategy $t\mapsto \Omega(t)$ depends
on the shape of the set $V$, and on the rate $M$ at which the contamination can be removed.
To state a precise result in this direction,  two geometric invariants must be introduced.

\begin{figure}[ht]
\centerline{\hbox{\includegraphics[width=15cm]{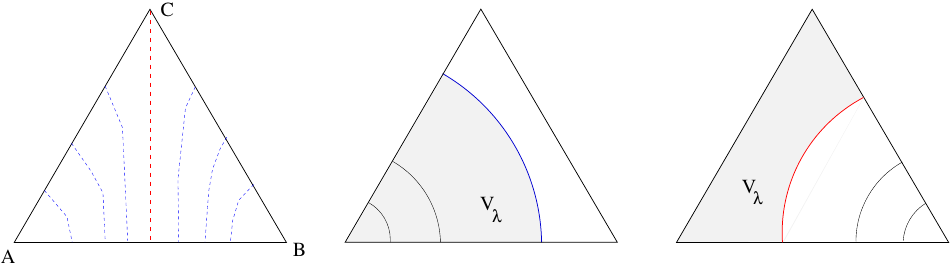}}}
\caption{\small The two invariants  (\ref{ko}) and (\ref{KO}) in the case of an equilateral triangle with unit side. Left: one of the level sets of the function $\phi$  in (\ref{KO}) 
must go through the vertex $C$. Hence it will have length $\geq \sqrt 3/2$.
Center: for any $\lambda\in [0, 1/2]$ we can cut a sector $V_\lambda$ with area 
$\lambda m_2(V)$, so that its boundary is an arc with length $~\leq~\sqrt{3\sqrt 3/4\pi}$.
Right: for  $\lambda\in [1/2, 1]$ we  can take $V_\lambda$ to be the complement 
of a sector with the same property.
}
\label{f:sc28}
\end{figure}

\begin{figure}[ht]
\centerline{\hbox{\includegraphics[width=13cm]{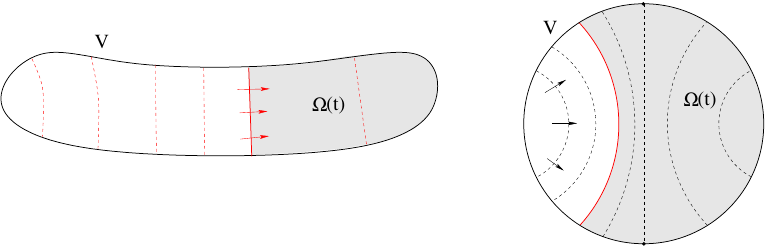}}}
\caption{\small  Left: if $M>K(V)$, an eradication strategy can be constructed in the form (\ref{erads}).
  Right: if $V$ is a disc, then the 
constants in (\ref{ko})-(\ref{KO}) are equal, and coincide with the diameter~$\delta$.
In this case, an eradication strategy exists if and only if $M> K(V)=\kappa(V)=\delta$.}
\label{f:sc42}
\end{figure}

\begi\item[(i)] First,  given $\lambda\in [0,1]$, among all subsets
$V_\lambda\subseteq V$ with area $m_2(V_\lambda)= \lambda \,m_2(V)$ we  
try to minimize the length of the  relative boundary
$\partial V_\lambda\cup V$.   This leads to the function
\bel{kl} \kappa(V,\lambda)~\doteq~\inf\Big\{ m_{1} \bigl( \partial V_\lambda\cap 
 V\bigr)\,;~~~V_\lambda\subseteq V, ~~m_2(V_\lambda)= \lambda \,m_2(V)\Big\}.\eeq
 Taking the supremum over all $\lambda\in [0,1]$, we define the constant
\bel{ko}\kappa(V)~\doteq~\sup_{\lambda\in [0,1]} \kappa(V,\lambda)
.\eeq
\v
\item[(ii)] Next, we  slice the set $V$  in terms of a continuous map $\phi:V\mapsto [0,1]$.
Here the slices are the pre-images $\phi^{-1}(\lambda)$, $\lambda\in [0,1]$.
We seek a function  $\phi$ such that  that the maximum size of the slices is as small as possible.
This yields the new constant (see \cite{Gromov, Guth})
\bel{KO}K(V)~\doteq~\inf_{\phi\,:\,V~\mapsto~ [0,1]} ~\left(
\sup_{\lambda\in [0,1]}  m_{1}\bigl( \phi^{-1}(\lambda)\bigr)\right).
\eeq
\endi
One always has $\kappa(V)\leq K(V)$. For some  sets,  a 
strict inequality   holds.    
\begin{example} {\rm If $V$ is a disc with diameter $\delta>0$, then $\kappa(V)= K(V)=\delta$.
In this case, one can prove that the contamination can be eradicated if and only if $M>\delta$ (see Fig.~\ref{f:sc42}, right).}
\end{example}
\begin{example} {\rm Let $V$ be an equilateral triangle whose side has unit length (Fig.~\ref{f:sc28}).
Then 
\bel{Ktri} K(V)~=~ {\sqrt 3\over 2}~\approx~0.866\ldots\eeq
coincides with the height of the triangle.  On the other hand, for 
any $\lambda\in [0,1]$ we have
\bel{ktri}\kappa(V, \lambda) ~\leq~\kappa\Big(V,{1\over 2}\Big)~=~     \sqrt {3\sqrt 3\over 4\pi}~=~\kappa(V)~\approx~0.643\ldots\eeq
}
\end{example}

The following theorem states that, for a general set $V$,
the contamination can be eradicated
if $M>K(V)$, and it cannot be eradicated if $M<   \kappa(V)$.

\begin{theorem}\label{t:61} Let $V\subset\R^2$ be a bounded open set with finite perimeter. 
Consider the following statements:
\begi
\item[(i)] There exists a continuous map $\vp: V\mapsto [0,1]$ whose level sets satisfy
\bel{slice}m_{1}\Big(\{ x\in V\,;~~\vp(x)=s\}\Big)~<~M\qquad\forall s\in [0,1].\eeq
\item[(ii)]  The constrained null controllability problem on $V$ has a solution.
\item[(iii)] For every $\lambda\in [0, 1]$ there exists a subset $V_\lambda\subset V$
such that 
\bel{cut}m_2(V_\lambda)~=~\lambda m_2(V),\qquad m_1(V\cap
\partial V_\lambda)~\leq~M.\eeq
\endi
Then we have the implications~ (i) $\implies$ (ii) $\implies$ (iii).
\end{theorem}
For a proof we refer to the forthcoming paper \cite{BMS}. 
The main ideas are illustrated in
 Figure~\ref{f:sc42}. 
 
 If $K(V)<M$, there exists a continuous map $\vp:V\mapsto [0,1]$ satisfying (\ref{slice}).
An eradication strategy can then be constructed in the form (see Fig.~\ref{f:sc42}, left)
\bel{erads}\Omega(t)~=~\bigl\{ x\in V\,;~~\vp(x)\geq \lambda(t)\bigr\},\eeq
for a suitable, strictly increasing map $\lambda:[0,T]\mapsto [0,1]$.
Indeed, at each time $t$ the relative boundary $\partial \Omega(t)\cap V$
has measure strictly less than $M$.   Hence  by 
\bel{dar2} \bega{rl}\ds{d\over dt} m_2\bigl(\Omega(t)\bigr)&\ds
=~\int_{\partial \Omega(t)\cap V} - \beta(t,x) \, d\sigma~=~
\int_{\partial \Omega(t)\cap V} \Big[1 - E\bigl(\beta(t,x)\bigr)\Big] \, d\sigma\\[4mm]
\ds &=~m_1\bigl(\partial \Omega(t)\cap V\bigr)-M~<~0,\enda\eeq
showing that
the area of $\Omega(t)$ is strictly decreasing.  In finite time, this area will shrink to zero.

On the other hand, if $M<\kappa(V)$, there exists $\lambda^*\in \,]0,1[\,$ such that, 
for every set $U\subseteq V$ having measure
$m_2(U) = \lambda^* m_2(V)$, the relative boundary is large:
$$m_1(\partial U\cap V)~>~M.$$
For any eradication strategy $t\mapsto\Omega(t)$, by continuity there will be a time $\tau\in [0,T]$
such that
$m_2\bigl(\Omega(\tau)\bigr)=  \lambda^* m_2(V)$.  
By (\ref{dar2}) this implies
$${d\over dt} m_2\bigl(\Omega(t)\bigr)\bigg|_{t=\tau}~=~m_1\bigl(\partial \Omega(\tau)\cap V\bigr)-M~>~0.$$
Hence the area of the sets $\Omega(t)$ can never shrink below $\lambda^* m_2(V)$, a contradiction.

\v

\section{Minimum time problems}
\label{sec:7}
\setcounter{equation}{0}
In the setting (\ref{EP}) considered in the previous section, the optimal solutions to the 
minimum time problem {\bf (MTP)} have been characterized in two basic cases.

\begin{figure}[ht]
\centerline{\hbox{\includegraphics[width=7cm]{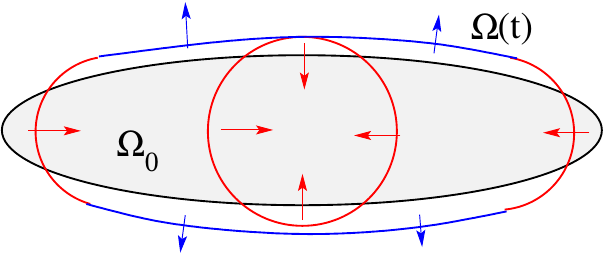}}\qquad\hbox{\includegraphics[width=6cm]{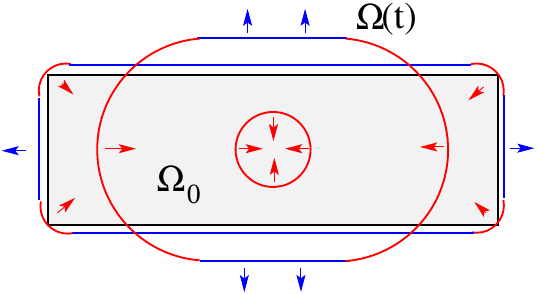}}}
\caption{\small  Optimal eradication strategies in the case where the initial set $\Omega_0$ is an ellipse, or a rectangle. At each time $t\in \,]0,T[\,$, the control is active along the arcs of circumferences, all with the same radius, where the curvature is maximum.  On the remaining portion of the boundary the control is not active and  the set expands.  In both examples, the boundary has $\C^{1,1}$ regularity, for $0<t<T$.  }
\label{f:sm68}
\end{figure}

The first result, proved in \cite{BBC}, deals with the eradication problem on the entire plane 
$\R^2$ (i.e., without geographical constraints), assuming that the
initial set $\Omega_0$ is convex  (see Fig.~\ref{f:sm68}).

\begin{theorem}\label{t:71} Let $\Omega_0\subset\R^2$ be a bounded open convex set, and assume that
$t\mapsto\Omega(t)$, $t\in [0,T]$, is a solution to the minimum time problem {\bf (MTP)}. 
Then the following holds.
\begi
\item[(i)] Every set $\Omega(t)$ is convex. For $0<t<T$, the boundary $\partial \Omega(t)$ has $\C^{1,1}$
regularity. Namely: the normal vector varies in a Lipschitz continuous way. In particular, the boundary has bounded curvature.
\item[(ii)] For every $t\in \,]0,T[\,$, the boundary control is active (i.e.: $\beta(t,x)>-1$)
on a finite or countable family of arcs
of circumferences, all with the same radius, where the boundary has maximum curvature.   
On the remaining part of the boundary the control is not active,
(i.e.: $\beta(t,x)=-1$) and the set expands with unit speed in the outer normal direction.
\endi
\end{theorem} 

The second  result yields a sufficient condition for the optimality of an eradication strategy
in the presence of a geographical constraint: $\Omega(t)\subseteq V$, for some bounded 
open set $V\subset\R^2$.    As before,  we assume that the effort 
$E$ is the function in (\ref{EP}).

\begin{theorem} \label{t:72} 
Let $t\mapsto \Omega(t)\subseteq V$ be an admissible strategy for the constrained eradication problem 
{\bf (CEP)}, with $\Omega(t)\subseteq \Omega(s)$ for all $0\leq s\leq t\leq 1$. 
Assume that, for all $0<t<T$, the following conditions hold.
\bel{a2}{d\over dt} m_2(\Omega(t))~=~m_1\bigl(V\cap \partial \Omega(t)\bigr)-M,\eeq
\bel{opc}m_1 \bigl(V\cap \partial \Omega(t)\bigr)~=~\min\Big\{ m_1 (V\cap \partial W)\,;
\qquad
W\subseteq V,\quad m_2(W) = m_2\bigl(\Omega(t)\bigr)\Big\}.\eeq
Then the strategy is optimal for the minimum time problem.
\end{theorem}
\v
\begin{figure}[ht]
\centerline{\hbox{\includegraphics[width=6cm]{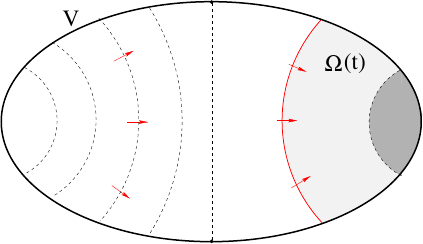}}}
\caption{\small An optimal solution for the constrained eradication problem in minimum time, in the case where
$V$ is an ellipse. For each $t\in \,]0,T[\,$, the
boundary $\partial \Omega(t)$ is an arc of circumference, crossing the boundary 
$\partial V$ perpendicularly.}
\label{f:sm4}
\end{figure}

\begin{remark} {\rm
Notice that the minimality condition in (\ref{opc}) implies that, for every 
$t\in [0,T]$, the relative boundary of $\Omega(t)$ has minimum length, compared with all other 
subsets of $V$ having the same area.  Therefore $\Omega(t)$  provides a solution to the classical Dido's problem
\cite{M}.   By the necessary conditions for optimality, the boundary $\partial \Omega(t)$ must be an arc of circumference, perpendicular to the boundary of $V$ at both endpoints
(see Fig.~\ref{f:sm4}).}
\end{remark}

\begin{figure}[ht]
\centerline{\hbox{\includegraphics[width=15cm]{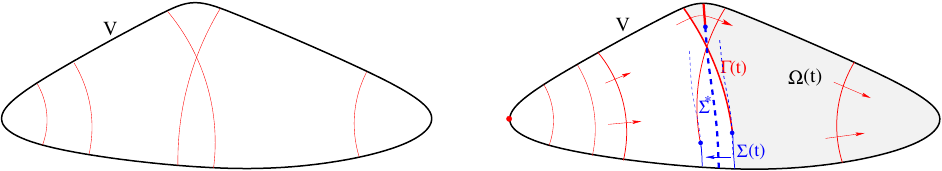}}}
\caption{\small  Left: a domain $V$ to which the sufficient conditions in Theorem~\ref{t:72} do not apply.
Indeed, arcs of circumferences which intersect the boundary $\partial V$ perpendicularly at both endpoints,
sometimes cross each other. Right: sketch of a (conjectured) optimal eradication strategy, constructed according 
to Steps 1--3 below.    At an intermediate time $t$, the control
is active along the arc of circumference $\Gamma(t)$, 
and not active along the remaining portion of boundary $\Sigma(t)$. }
\label{f:sm75}
\end{figure}

Unfortunately, there are many domains $V\subset\R^2$ for which the conditions in Theorem~\ref{t:72} 
are not satisfied
by any map $t\mapsto\Omega(t)\subseteq V$.    
For example, think of a triangle. Another example is shown in
Fig.~\ref{f:sm75}, left.  This leads to our third main open problem:
\v
{\bf Problem 3.} {\it Given an open connected domain $V\subset\R^2$ with bounded perimeter,
describe the optimal strategy $t\mapsto \Omega(t)\subseteq V$ for {\bf (CEP)},  achieving the eradication 
in minimum time.
In particular, it would be of interest to study the cases where  (i) $V$ is a polygon, or (ii) $V$ is a convex set with smooth boundary.}
\v

Assuming that the boundaries $\partial \Omega(t)$ are sufficiently regular, some insight can be gained
by the analysis of optimality conditions.
When the effort function is $E(\beta)= \max\{ 1+\beta, 0\}$  and $\kappa_1=0$, the adjoint equations (\ref{Ydt}) 
take the simple form
\bel{aey}
Y_t(t,\xi)~=~-\omega(t,\xi) Y(t,\xi),\eeq
where $\omega(t,\xi)$ is the curvature of the boundary $\partial \Omega(t)$ at the point $x(t,\xi)$.
According to Theorem~6.1 in \cite{BCS2}, if $t\mapsto\Omega(t)$ is a strategy that
minimizes the terminal area $m_2\bigl(\Omega(T)\bigr)$, then there exists   a scalar function
 $t\mapsto \lambda(t)> 0$
such that  the normal velocity $\beta=\beta(t,\xi)$ satisfies
\bel{max16}
\lambda(t) E\big(\beta(t,\xi)\bigr) -Y(t,\xi) \beta(t,\xi)  ~= ~\min_{\beta\geq -1}~
\Big\{\lambda(t) E(\beta)-Y(t,\xi) \beta  \Big\}.\eeq
Notice that, if the minimum is achieved by some value $\beta(t,\xi)>-1$, then necessarily 
$Y(t,\xi)=\lambda(t)$ is independent of $\xi$.  

In view of (\ref{aey}), the requirement  that
$Y(t,\xi) = \lambda(t)$ is a function of time alone
 implies the constant curvature condition:
\begi
\item[{\bf (CC)}]
{\it At any time $t\in [0,T]$ the curvature $\omega(t,\cdot)$ must  be constant along the portion of the boundary 
where the control is active, i.e.: where $\beta(t,\xi)>-1$.}
\endi
The picture that emerges from the study of optimality conditions is the following:
\begi
\item At a.e.~time $t\in [0,T]$, the relative boundary of the set $\Omega(t)$ can be decomposed as
$$V\cap \partial \Omega(t)~=~\Sigma(t) \cap \Gamma(t),$$
where $\Sigma(t)$ is the portion where no control effort is made: $\beta=-1$, and the set thus expands at unit rate,
while $\Gamma(t)$ is the portion where the control is active: $\beta>-1$.
\item $\Gamma(t)$ is the union of arcs of circumferences, all with the same radius $r(t)$.
\item The arcs in $\Gamma(t)$ should join $\Sigma(t)$ tangentially, and cross the boundary $\partial V$ perpendicularly.
\endi

The above considerations suggest a possible approach to determine the optimal motion $t\mapsto \Omega(t)$
(see Fig.\ref{f:sm75}, right).

{\bf Step 1.}  Construct the uncontrolled portion $\Sigma^*= \Sigma(t^*)$ of the boundary at a particular time $t^*$ 
where its extension  is maximal. 

{\bf Step 2.} Construct the uncontrolled portions $\Sigma(t)$, for $t$ in a neighborhood of  $t^*$, 
observing that they satisfy
 $$\Sigma(t)~\subset~\Big\{ x\in V\,;~~d\bigl(x, \Sigma(t^*)\bigr)~=~|t-t^*|\Big\}.$$

{\bf Step 3.}  Construct the remaining portion of the boundary $\Gamma(t)$,  
observing that it consists of arcs of circumferences all with the same radius, 
tangent to $\Sigma(t)$ and perpendicular to $\partial V$.
The position of these circumferences should be uniquely determined by the constraint on the total effort
(\ref{a2}).
\v
In this approach, the key issue (still poorly understood) 
is how to determine the ``maximal free boundary" $\Sigma^*$, by solving
a suitable system of ODEs together with boundary conditions.

\begin{figure}[ht]
\centerline{\hbox{\includegraphics[width=15cm]{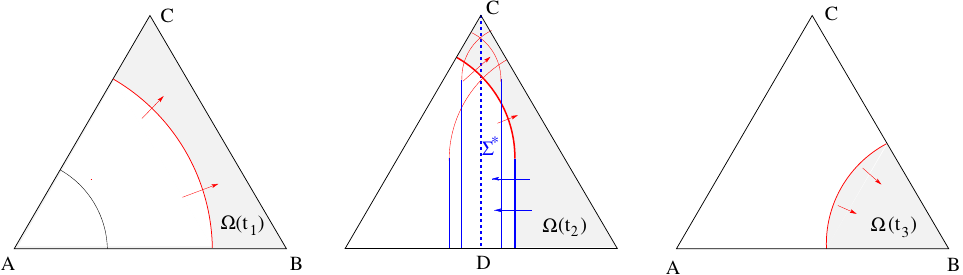}}}
\caption{\small  The optimal eradication strategy, when $V$ is the isosceles triangle $ABC$.
At the intermediate time $t_2$, the relative boundary of the set $\Omega(t_2)$ is the union of a vertical segment where the 
control is not active: $\beta=-1$, and an arc of circumference where the control is active: $\beta>-1$. 
Here the maximal free boundary $\Sigma^*$ is the vertical segment $CD$. }
\label{f:sc26}
\end{figure}

Figure~\ref{f:sc26} shows the construction of an optimal strategy when $V$ is an isosceles triangle.
In this case, by symmetry, it is easy to guess that $\Sigma^*$ should be
the vertical segment through the vertex $C$.

\begin{figure}[ht]
\centerline{\hbox{\includegraphics[width=9cm]{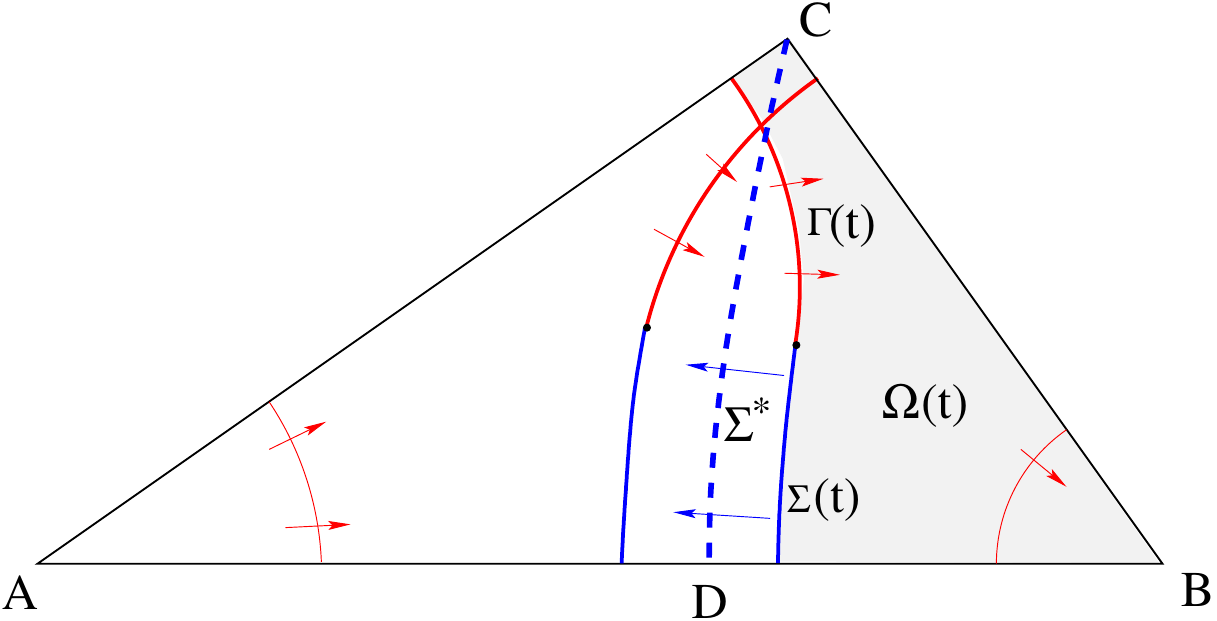}}}
\caption{\small  Sketch of the optimal eradication strategy, for a general triangle $ABC$. Here  the 
maximal free boundary is the dotted curve $\Sigma^*$. It is tangent to the bisectrix of the angle at $C$, and perpendicular to the segment $AB$ at the point $D$.
At an intermediate time $t$, the relative boundary of the shaded region $\Omega(t)$ is
the union of (i) an arc $\Sigma(t)$ parallel to $\Sigma^*$, where the 
control is not active: $\beta=-1$, and (ii) an arc of circumference $\Gamma(t)$ where the control is active: $\beta>-1$. 
 }
\label{f:sm74}
\end{figure}

A similar construction was carried out  in \cite{BMS}  for a general triangle (see Fig.~\ref{f:sm74}).
In this case the maximal free boundary $\Sigma^*$ is determined by solving a second
order, implicit, singular system of ODEs.

\v

Before closing, we wish to point out a class of somewhat simpler geometrical problems, whose solutions
 may provide some additional insight.

Starting with the minimum time eradication problem and taking the limit 
$M\to \infty$, after a time rescaling $t' = Mt$ one formally obtains

\begi
\item[{\bf (OSP)}] {\bf (Optimal slicing problem).} {\it Let $V\subset\R^2$ be a compact connected set with 
$m_2(V)= A >0$.
Construct a family of sets $t\mapsto \Omega(t)\subseteq V$ whose relative boundaries minimize the
integral
\bel{idom}
\int_0^A m_1\bigl(V\cap \partial \Omega(t)\bigr)\, dt,\eeq
subject to
\bel{m2o}m_2\bigl(\Omega(t)\bigr) ~=~t,\qquad \qquad\forall t\in [0,A],\eeq
\bel{cresc}\Omega(t_1)\subset \Omega(t_2),\qquad\qquad\hbox{for}~~ 0\leq t_1<t_2\leq A.\eeq
}
\endi

\begin{figure}[ht]
\centerline{\hbox{\includegraphics[width=11cm]{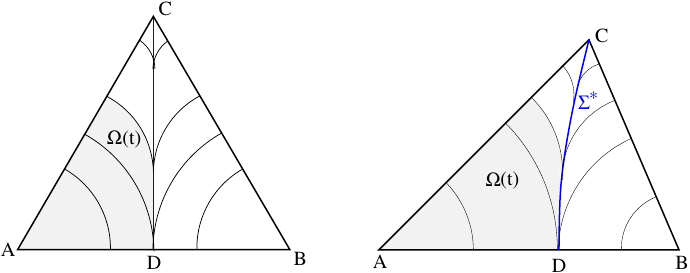}}}
\caption{\small  A sketch of the optimal slicing for an isosceles triangle (left) and a general triangle (right).
Here $\Sigma^*$ is a curve tangent to the bisectrix of the angle at $C$, and perpendicular to the segment $AB$ at $D$. The boundary of $\Omega(t)$ contains a portion of $\Sigma^*$, and an arc of circumference which is
tangent to $\Sigma^*$ at one endpoint, and perpendicular to one of the sides of the triangle at the other endpoint.
}
\label{f:sm72}
\end{figure}

One may regard {\bf (OSP)} as a ``time dependent Dido's problem", where we seek to minimize the 
average length of the boundaries of the sets $\Omega(t)$.    If the set $V$ has a shape to which Theorem~\ref{t:72} applies, then  (after a rescaling of time) the same sets $\Omega(t)$  in (\ref{a2})-(\ref{opc})
provide a solution to the optimal slicing problem.  
In the case of triangles, solutions to {\bf (OSP)} are sketched in Fig.~\ref{f:sm72}.
It remains an open problem to understand the optimal slicing of a set such as the one in Fig.~\ref{f:sm75}.

\end{document}